\documentclass[10pt,a4paper]{article}

\usepackage[utf8]{inputenc}
\usepackage[english]{babel}
\usepackage{amsmath}  
\usepackage{amsfonts} 
\usepackage{amsthm}    
\usepackage{amssymb}
\usepackage{graphicx, overpic} 
\usepackage{paralist} 
\usepackage{geometry} 
\usepackage{xcolor, comment}
\usepackage{bbm}
\usepackage{hyperref}

\author{Christian Bick${}^{1,2,3,4}$ and Tobias Böhle${}^{1,2,5}$ and Christian Kuehn${}^{5,6,7}$}
\date{{\small
${}^1$Department of Mathematics, Vrije Universiteit Amsterdam, De Boelelaan 1111, Amsterdam, the Netherlands.\\
${}^2$Institute for Advanced Study, Technical University of Munich, Lichtenbergstr. 2, 85748 Garching, Germany.\\
${}^3$Department of Mathematics, University of Exeter, Exeter EX4 4QF, United Kingdom.\\
${}^4$Mathematical Institute, University of Oxford, Oxford OX2 6GG, United Kingdom.\\
${}^5$Technical University of Munich, School of Computation Information and Technology, Department of Mathematics, Boltzmannstr. 3, 85748 Garching, Germany.\\
${}^6$Complexity Science Hub Vienna, Josefst\"adter Str.~39, 1080 Vienna, Austria.\\
${}^7$Munich Data Science Institute, Walther-von-Dyck-Str. 10, 85748 Garching, Germany.\\}
\medskip
\today}
\title{Higher-Order Network Interactions through Phase Reduction for Oscillators with Phase-Dependent Amplitude}

\newcommand{\R}{\mathbb{R}} 
\newcommand{\C}{\mathbb{C}} 
\newcommand{\Z}{\mathbb{Z}} 
\newcommand{\N}{\mathbb{N}}
\renewcommand{\S}{\mathbb{S}}
\renewcommand{\top}{\mathsf{T}}

\newcommand{\abs}[1]{\left\lvert #1 \right\rvert}

\renewcommand{\d}{\mathrm d}

\theoremstyle{definition}

\theoremstyle{plain} 	

\theoremstyle{remark}

\theoremstyle{plain}

\numberwithin{equation}{section}

\begin{document}

\maketitle

\hrule
\paragraph{Abstract.}
Coupled oscillator networks provide mathematical models for interacting periodic processes. 
If the coupling is weak, phase reduction---the reduction of the dynamics onto an invariant torus---captures the emergence of collective dynamical phenomena, such as synchronization.
While a first-order approximation of the dynamics on the torus may be appropriate in some situations, higher-order phase reductions become necessary, for example, when the coupling strength increases.
However, these are generally hard to compute and thus they have only been derived in special cases: This includes globally coupled Stuart--Landau oscillators, where the limit cycle of the uncoupled nonlinear oscillator is circular as the amplitude is independent of the phase.
We go beyond this restriction and derive second-order phase reductions for coupled oscillators for arbitrary networks of coupled nonlinear oscillators with phase-dependent amplitude, a scenario more reminiscent of real-world oscillations. 
We analyze how the deformation of the limit cycle affects the stability of important dynamical states, such as full synchrony and splay states.
By identifying higher-order phase interaction terms with hyperedges of a hypergraph, we obtain natural classes of coupled phase oscillator dynamics on hypergraphs that adequately capture the dynamics of coupled limit cycle oscillators.
\vspace{-3mm}
 \vspace{-3mm}
\vspace{4mm}
\hrule

\section{Introduction}

Collective behavior of oscillatory systems, such as synchronization, is ubiquitous in many real-world dynamical systems. 
Examples range from biological systems, such as the continuous beat of a healthy heart, regularly spiking neurons, or flashing fireflies to large scale technological systems, such as power grids, or the periodic motion of planets~\cite{Pikovsky2001, Buck1968, Golomb2001}.
If the coupling is sufficiently weak, phase reductions provide a useful tool to study these oscillator networks~\cite{Pietras2019, Monga2019, Nakao2016} and have found application to elucidate collective dynamics for example in neuroscience~\cite{Ashwin2016}.
Specifically, if one assumes that each oscillatory unit has a stable limit cycle without coupling, then the whole system still settles to an invariant torus for sufficiently small coupling strengths.
A phase-reduced system then describes the dynamics on this invariant torus and it is usually derived as an expansion in the coupling strength. 
Typically, one considers a first-order truncation ignoring any higher-order terms. However, this is insufficient to describe the dynamics of the full/unreduced system when the first-order truncation undergoes a bifurcation or when the coupling is stronger.

To accurately describe the dynamics of the full system in these cases, one has to resort to higher-order phase reductions. 
Recently progress has been made to compute such phase reductions: 
Explicit computations show how nonpairwise phase interactions enter the phase-reduced equations once one goes to second or higher orders~\cite{Leon2019, Gengel2021}. 
In general, however, computing higher-order phase reductions is not straightforward and the focus has been on simple oscillator models with additive interactions.
For example, the computations in~\cite{Leon2019} make explicit use of the rotational symmetry of Stuart--Landau oscillators that imply that the (unperturbed) limit cycle is the unit circle; the resulting phase equations reflect the symmetry properties.
Such an assumption of rotational symmetry is rarely satisfied for general oscillator systems.
Indeed, while a limit cycle may be approximately circular for oscillations emanating from a Hopf bifurcation point, oscillations further away from the bifurcation point will generically have a phase-dependent amplitude.


Here, we derive higher-order phase reductions for systems in which the limit cycle has phase-dependent amplitude. More specifically, we generalize recent approaches for coupled Stuart--Landau oscillators on a graph (i.e., coupling is additive) to oscillators subject to a perturbation that breaks the rotational symmetry of the limit cycle.
We derive phase equations by expanding in terms of both the coupling strength between oscillators as well as the size of the symmetry-breaking perturbation. 
We then analyze how these higher-order interaction terms affect the stability full synchrony---all oscillators are at the same state---and the splay configuration, in which the phases of the oscillators are uniformly spread out.
Our approach allows to compute phase reductions not only for symmetric all-to-all coupled networks, but also coupled oscillators on arbitrary graphs.
Thus, for coupled oscillators on a given graph with additive interactions, we obtain a parameterized family of effective phase dynamics which include higher-order nonpairwise phase interactions that depend nonlinearly on three or more oscillator phases.

Our results also provide a tool to construct phase dynamics on hypergraphs that are a meaningful approximation of systems of nonlinearly coupled oscillators.
Indeed, nonpairwise phase coupling terms between more than two oscillator phases have been associated with phase oscillator dynamics on hypergraphs~\cite{Battiston2020,Bick2021a}; these can arise for example through higher-order phase reductions of additively coupled systems (e.g.,~\cite{Leon2019}) or first-order phase reductions of oscillators with generic nonlinear coupling~\cite{Ashwin2016}.
So far, many phase oscillator networks with higher-order interactions that have been considered were ad-hoc, for example, by generalizing the Kuramoto model to hypergraphs (see, e.g.,~\cite{Skardal2019, Bick2022}).
By contrast, our results provide a natural family of hypergraphs together with phase interaction functions that describe the synchronization behavior of an (unreduced) nonlinear oscillator network.
This family is parameterized in terms of the underlying coupling graph as well as the system parameters.

This paper is organized as follows: 
In Section~\ref{sec:prelim} we recall the main points of a previous work~\cite{Leon2019} considering higher-order phase reductions for globally coupled Stuart--Landau oscillators. 
Then, in Section~\ref{sec:meat} we study a system whose limit cycle can be obtained from perturbing the circular limit cycle from a Stuart--Landau oscillator. 
We derive phase reductions as an expansion in the coupling strength up to second order and the parameter that controls the deformation of the limit cycle. 
Next, in Section~\ref{sec:dynamics} we numerically analyze how the deformation of the limit cycle affects the stability of synchronized and splay states. 
We investigate how accurately first and second order phase reductions reproduce these stability properties. 
Finally, Section~\ref{sec:conclusion} contains a discussion and some concluding remarks.

\section{Phase Reductions for Stuart--Landau Oscillators}\label{sec:prelim}

In this section, we recall the main aspects of how to derive higher-order phase reductions for coupled Stuart--Landau oscillators from~\cite{Leon2019}. We highlight the main assumptions that are made to derive these reductions.

First, let us consider a single complex Stuart--Landau oscillator with state $A=A(t)\in \C$, that evolves according to
\begin{align}\label{eq:LP_single}
	\dot A :=  \frac{\d}{\d t}A = A - (1+ic_2)\abs{A}^2 A,
\end{align}
where $c_2\in\R$ is a parameter. The right-hand side of~\eqref{eq:LP_single} is equivariant with respect to the continuous group~$\S:= \R/(2\pi \Z)$, which acts on~$\C$ by shifting an oscillator by a given phase. Due to this symmetry, it makes sense to introduce polar coordinates $A=re^{i\phi}$ with $r\ge 0$ and $\phi\in \S$. Then, the system has a stable limit cycle at $r = 1$ and the dynamics on this limit cycle can be described by just the phase $\phi$. In fact, on the limit cycle, the phase $\phi$ increases with constant speed $-c_2$, such that $\phi(t) = \phi(0) - c_2t$. To understand the dynamics off the limit cycle, we note that every point in the basin of attraction of this attractive limit cycle has an asymptotic phase, which is the phase of an initial condition condition of a trajectory on the limit cycle that the point converges to. The set of points with the same asymptotic phase are called isochrons~\cite{Guckenheimer1975, Langfield2014}. To define them, we introduce the notation $\Phi^t A_0$ for the solution of~\eqref{eq:LP_single} with $A(0) = A_0$ at time $t$. If $A_0$ is on the limit cycle, i.e., $\abs{A_0}=1$, and $\phi(t) = \theta- c_2t$ then
\begin{align*}
	\mathcal I(\theta) := \{ \hat A_0 \in \C: \lim_{t\to\infty} \Phi^t \hat A_0 - e^{i(\theta - c_2t)} = 0\}
\end{align*}
is the isochron with asymptotic phase $\theta$. Upon variation of $\theta$, the isochrons foliate the basin of attraction of the limit cycle. For~\eqref{eq:LP_single}, they can explicitly be calculated~\cite{Leon2019} to be
\begin{align*}
\mathcal I(\theta) =  \{ A = re^{i\phi}: \theta = \phi-c_2\ln r\}.
\end{align*}
As one can see from this formula these isochrons are symmetric in the sense that one isochron can be obtained from another by shifting it by a constant phase. This property also follows directly from the~$\S$ symmetry of~\eqref{eq:LP_single}.

Having studied a single Stuart--Landau oscillator, we now consider $N\in \N$ coupled Stuart--Landau oscillators described by complex variables $A_k$, $k=1,\dots,N$. When the coupling is as in~\cite{Leon2019}, they satisfy
\begin{align}\label{eq:LP_CGL}
\dot A_k = A_k - (1+i c_2)\abs{A_k}^2A_k + \epsilon (1+ic_1)(\bar A-A_k), \quad k=1,\dots,N.
\end{align}
Here, $c_1,c_2$ are two real parameters, $\epsilon\ge 0$ relates to the coupling strength and $\bar A = \frac 1N \sum_{j=1}^N A_j$. Now one changes to polar coordinates $A_k = r_k e^{i\phi_k}$, with $r_k\ge 0$ and $\phi_k\in \S:= \R/(2\pi \Z)$. After conducting an additional nonlinear transformation $\theta = \phi-c_2 \ln(r)$ to straighten the isochrons, the authors of~\cite{Leon2019} arrive at the system
\begin{subequations}
\label{eq:LP_system}
\begin{align}
\label{eq:LP_system_radius}
\dot r_k &= f(r_k) + \epsilon g_k(r,\theta),\\
\label{eq:LP_system_phase}
\dot \theta_k &= \omega + \epsilon h_k(r,\theta),
\end{align}
\end{subequations}
where $\omega\in \R$ is a parameter that depends on $c_2$, $f(r) = r(1-r^2)$ and the functions $g_k$ and $h_k$ are given by
\begin{align*}
g_k(r,\theta) &= -r_k + \frac 1N \sum_{j=1}^N \left\{r_j\left[\cos\left( \theta_j-\theta_k + c_2 \ln \frac{r_j}{r_k}\right) - c_1 \sin\left( \theta_j-\theta_k + c_2 \ln \frac{r_j}{r_k}\right)\right]\right\},\\
h_k(r,\theta) &= c_2-c_1 + \frac{1}{Nr_k}\sum_{j=1}^N\Bigg\{ r_j\bigg[(c_1-c_2)\cos\left(\theta_j-\theta_k+c_2\ln\frac{r_j}{r_k}\right)\\
&\qquad \qquad \qquad \qquad \qquad + (1+c_1c_2) \sin\left( \theta_j-\theta_k + c_2 \ln\frac{r_j}{r_k}\right)\bigg]\Bigg\}.
\end{align*}
The coupling in~\eqref{eq:LP_CGL} respects the $\S$ symmetry of a Stuart--Landau oscillator such that~\eqref{eq:LP_CGL} again possesses an $\S$ symmetry group that acts on $\C^N$ by shifting all oscillators $A_1,\dots,A_N$ by the same phase. Since the transformation that straightens the isochrons does not break this symmetry, the system~\eqref{eq:LP_system} inherits the same symmetry group. This fact can also be observed by directly looking at the structure of the functions $g_k, h_k$ and $f$. Since they only depend on phase differences, they are invariant when all oscillators are shifted by the same phase. Therefore, we can change into a co-rotating coordinate frame
\begin{align}\label{eq:LP_S1_trafo}
	\theta \mapsto \theta + \omega t,
\end{align}
and thereby set $\omega=0$ without loss of generality.

Next, we derive first and second order phase reductions for the system~\eqref{eq:LP_system}. In absence of the coupling, i.e., $\epsilon= 0$, each oscillator in~\eqref{eq:LP_CGL} and~\eqref{eq:LP_system} has a stable limit cycle at $\abs{A_k} = 1$ or $r_k=1$, respectively. Therefore, the limiting dynamics of the whole system takes place on the $N$-dimensional torus that is described by $r_k\equiv 1$ for all $k=1,\dots,N$. When slightly increasing $\epsilon$ this torus persists but it gets perturbed. The radii of this invariant torus are then functions of the phases. In fact, they can be expanded in terms of $\epsilon$ such that
\begin{align}\label{eq:LP_r_expansion}
	r_k(\theta) = r_k^{(0)}(\theta) + \epsilon r_k^{(1)}(\theta) + \epsilon^2 r^{(2)}(\theta) + \mathcal O(\epsilon^3),
\end{align}
with $r_k^{(0)} \equiv 1$, see~\cite{Leon2019}. Inserting the ansatz~\eqref{eq:LP_r_expansion} into~\eqref{eq:LP_system_phase} as done in~\cite{Leon2019} yields
\begin{align}\label{eq:LP_phase_reduction}
	\dot \theta_k = \epsilon h_k(r^{(0)}, \theta) + \epsilon^2 \nabla_r h_k(r^{(0)}, \theta) \cdot r^{(1)}(\theta) + \mathcal O(\epsilon^3),
\end{align}
since $\omega = 0$. By truncating terms of order $\mathcal O(\epsilon^2)$ and higher orders, one can obtain a first-order phase reduction. To derive a second order phase reduction one also needs to know $r^{(1)}(\theta)$. This can be accomplished by inserting~\eqref{eq:LP_r_expansion} into~\eqref{eq:LP_system_radius} and collecting terms of order $\epsilon$, see~\cite{Leon2019}. One then ends up with
\begin{align}\label{eq:LP_eq1}
	\dot r_k^{(1)} = f'(r_k^{(0)}) r^{(1)}_k + g_k(r^{(0)},\theta).
\end{align}
Moreover, by the chain rule, one has
\begin{align*}
	\dot r_k = (\nabla_\theta r_k)\cdot \dot \theta &= (\nabla_\theta r_k)\cdot (\omega\mathbbm 1 + \epsilon h(r,\theta))\\
	&= \Big(\nabla_\theta r_k^{(0)} + \epsilon \nabla_\theta r_k^{(1)} + \mathcal O(\epsilon^2)\Big) \cdot \Big(\omega\mathbbm 1 + \epsilon h(r^{(0)} + \mathcal O(\epsilon),\theta)\Big),
\end{align*}
for the dynamics on the invariant torus, where $\mathbbm 1 = (1,\dots,1)^\top \in \R^N$. Now, it is crucial that $\omega$ can be set to $0$, because then collecting terms of order $\epsilon$ in this equation yields
\begin{align}\label{eq:LP_eq2}
	\dot r^{(1)}_k = (\nabla_\theta r_k^{(0)})\cdot h(r^{(0)}, \theta).
\end{align}
Because $\nabla_\theta r^{(0)}_k = \nabla_\theta 1 = 0$, combining~\eqref{eq:LP_eq1} and~\eqref{eq:LP_eq2}, the authors of~\cite{Leon2019} arrive at
\begin{align*}
r^{(1)}_k = -\frac{g_k(r^{(0)},\theta)}{f'(r^{(0)})} = \frac 12 g_k(r^{(0)}, \theta).
\end{align*}
Substituting that into~\eqref{eq:LP_phase_reduction} and truncating $\mathcal O(\epsilon^3)$ terms yields the second order phase reduction.

\section{Phase Reductions for Limit Cycles with Phase-Dependent Amplitude}\label{sec:meat}

In this section, we introduce a variation of Stuart--Landau oscillators where the limit cycle is not circular but has a phase dependent amplitude. We then derive first- and second-order phase reductions reductions for this class of oscillators subject to coupling as in the previous section.
Finally, we investigate how these phase reductions are affected by the parameter that determines the deviation of the shape of the limit cycle from a circle.

\begin{figure}[h]
    \centering
    \begin{overpic}[width = 0.5\textwidth]{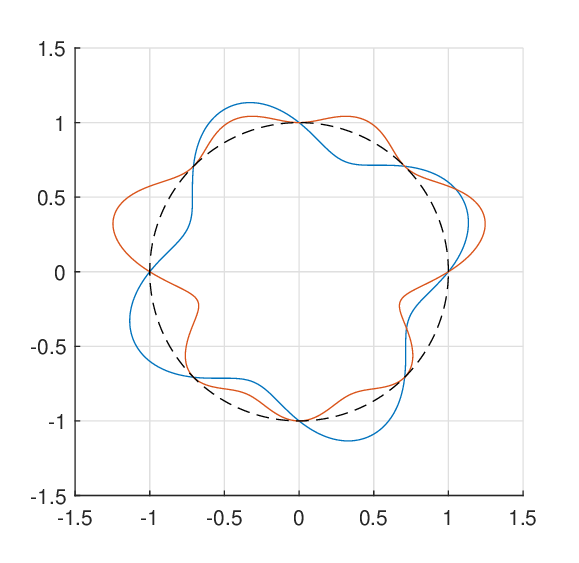}
    \put(45,2){$\operatorname{Re}(A)$}
    \put(-1,56){$\operatorname{Im}(A)$}
    \end{overpic}
    \caption{Different deformations of limit cycles captured by our system \eqref{eq:single_oscillator}. The blue curve represents the limit cycle when $g(\phi) = \sin(4\phi)$ and $\delta = 0.2$, while the red curve depicts the limit cycle when $g(\phi) = \sum_{i=1}^4 \sin((2i-1)\phi)$ and $\delta = 0.1$. These limit cycles are deviations from the unit circle (black dashed line).}
    \label{fig:deformed_limitcycle}
\end{figure}

Inspired by~\cite{Leon2019}, we start with a modified Stuart--Landau oscillator that can have a noncircular limit cycle of a given functional form. Specifically, for a parameter~$\delta\in\R$ with $|\delta|\ll 1$ and a given smooth $2\pi$-periodic function $g\colon \S\to\R$, the oscillator we consider has a limit cycle with phase-dependent amplitude $r = 1+\delta g(\phi)$.
This is the case if the state of oscillator $A = re^{i\phi}\in \C$ evolves according to
\begin{subequations}
\label{eq:single_oscillator}
\begin{align}
\label{eq:single_oscillator_r}
\dot r &= \delta g'(\phi) \omega \frac{r}{1+\delta g(\phi)} + m r^2 (r-1-\delta g(\phi)),\\
\dot \phi &=\omega,
\end{align}
\end{subequations}
where $\omega>0$ is the angular velocity of the oscillator and $m<0$ determines the rate of attraction to the limit cycle.
The shape of the limit cycle is shown in Figure \ref{fig:deformed_limitcycle}; for $\delta=0$ the limit cycle is circular.
Due to the explicit embedding of the limit cycle, we choose the nonlinearity of the radial direction to be slightly different as in the Stuart--Landau oscillator~\eqref{eq:LP_single}.
As this primarily influences the rate of radial convergence, this should be not relevant for the phase reductions since the derivatives of the radial equation at $r=1$ coincide if $m=-2$.
Moreover, to keep the problem tractable, we focus on oscillations without radial dependency of the phase dynamics (corresponding to $c_2=0$ in~\eqref{eq:LP_single}).

Now, given an ensemble of $N$ oscillators $(A_k = r_ke^{i\phi_k})_{k=1,\dots, N}$, we assume a mean-field coupling, of the form
\begin{align}\label{eq:coupling}
\dot A_k = \mathcal F(A_k) + K e^{i\alpha} (\bar A-A_k),
\end{align}
where $\mathcal F(A_k)$ denotes the intrinsic dynamics of oscillator $k$ as described in~\eqref{eq:single_oscillator}, $K\in \R$ is the coupling strength, $\alpha\in\S$ is a parameter and $\bar A = \frac 1N \sum_{j=1}^N A_j$ is the average position as before. Rewritten in polar coordinates, this results in the system
\begin{align*}
	\dot r_k &= \delta g'(\phi_k) \omega \frac{r_k}{1+\delta g(\phi_k)} + m r_k^2 (r_k-1-\delta g(\phi)) + \frac{K}{N}\sum_{l=1}^N [r_l\cos(\phi_l-\phi_k + \alpha) - r_k\cos(\alpha)]\\
	\dot \phi_k &= \omega + \frac{K}{N r_k}\sum_{l=1}^N [r_l\sin(\phi_l-\phi_k+\alpha) - r_k\sin(\alpha)].
\end{align*}
After the transformation
\begin{align*}
	R_k = \frac{r_k}{1+\delta g(\phi_k)},
\end{align*}
to transform the phase-dependent limit cycle to a circle, we arrive by a direct, but slightly lengthy, calculation at the system
\begin{subequations}
\label{eq:main_system}
\begin{align}
\label{eq:main_system_R}
\dot R_k &= F(R_k, \phi_k) + KG_k(R,\phi),\\
\label{eq:main_system_phi}
\dot \phi_k&= \omega + K H_k(R, \phi),
\end{align}
\end{subequations}
with functions $F, G_k$ and $H_k$ defined by
\begin{align*}
	F(R_k, \phi_k) &= m R_k^2(R_k-1)(1+\delta g(\phi_k))^2,\\
	G_k(R,\phi) &= \frac{1}{N}\sum_{l=1}^N\Bigg[ R_l\frac{1+\delta g(\phi_l)}{1+\delta g(\phi_k)}\cos(\phi_l-\phi_k+\alpha) - R_k\cos(\alpha)\\
	&\quad - \delta g'(\phi_k)\left( R_l\frac{1+\delta g(\phi_l)}{(1+\delta g(\phi_k))^2}\sin(\phi_l - \phi_k + \alpha) - R_k \frac{\sin(\alpha)}{1+\delta g(\phi_k)}\right) \Bigg],\\
	H_k(R,\phi) &= \frac{1}{N}\sum_{l=1}^N \left[ \frac{R_l(1+\delta g(\phi_l))}{R_k (1+\delta g(\phi_k))}\sin(\phi_l-\phi_k+\alpha) - \sin(\alpha)\right].
\end{align*}
It is important to note, that for a general perturbation $F, G$ and~$H$ depend explicitly on the phase. 
Therefore, the system~\eqref{eq:main_system} does not have an $\S$~symmetry. 
Hence, we cannot change into a co-rotating coordinate system. That means, unlike in the system~\eqref{eq:LP_system}, we cannot assume $\omega=0$ without loss of generality. 
However, $\omega=0$ was a crucial assumption to derive formula~\eqref{eq:LP_eq2} in Section~\ref{sec:prelim}. 
Next, we show how to generalize the methods of~\cite{Leon2019} to derive higher-order phase reductions anyway.

\subsection{Phase Reductions of First-Order}\label{sec:first_order}

Because there is an additional parameter~$\delta$ in our system, we expand in both~$K$ and~$\delta$, i.e., we study asymptotics in two parameters~\cite{Kuehn2022}.
Regarding notation, if $W$ is a function or a scalar, we write $W^{(n,j)}$ for the contribution of order $K^n\delta^j$ to~$W$, i.e.,
\begin{align*}
	W = \sum_{n,j=0}^\infty K^n \delta^j W^{(n,j)}.
\end{align*}
Moreover, we write $W^{(n, \star)}$ for all contributions of order $K^n$, which includes all orders in $\delta$. Similarly, $W^{(\star,j)}$ includes all terms of order $\delta^j$. Consequently,
\begin{align*}
	W = \sum_{j=0}^\infty \delta^j W^{(\star, j)} = \sum_{n=0}^\infty K^n W^{(n,\star)}.
\end{align*}
In particular, if the quantity $W$ is independent of $K$, we have $W=W^{(0,\star)}$, but we use the notation $W=W^{(-,\star)}$ to highlight the independence of $K$. We use this notation to derive phase reductions of different approximation order in $K$ and $\delta$. To distinguish these phase reductions, we speak of an $(a,b)$-phase reduction when $a$ is the highest approximation order in $K$ and $b$ is the highest order in $\delta$. In particular, an $(a,b)$-phase reduction is given by
\begin{align*}
	\phi_k = \sum_{n=0}^a \sum_{j=0}^b K^n \delta^j P_k^{(n,j)}(\phi),
\end{align*}
where $P_k^{(n,j)}(\phi)$ denotes the contribution on the order $K^n\delta^j$. Explicit expressions for $P_k^{(n,j)}(\phi)$ will be derived below.

By the same reasons as illustrated in Section~\ref{sec:prelim}, the limiting dynamics of the system~\eqref{eq:main_system} takes place on an attractive invariant $N$-dimensional torus.
If $K=0$ this torus is described by $R_k\equiv 1$ for all $k=1,\dots,N$. If $\abs{K}$ is small, the torus persists and the radii of this torus can be expanded in terms of $K$ as
\begin{align}\label{eq:R_expansion}
	R_k(\phi) = R_k^{(0,\star)}(\phi) + K R_k^{(1,\star)}(\phi) + K^2 R_k^{(2,\star)}(\phi) + \mathcal O(K^3),
\end{align}
where $R_k^{(0,\star)}(\phi)\equiv 1$.
When inserting the expansion~\eqref{eq:R_expansion} into the system~\eqref{eq:main_system_phi} we obtain
\begin{align}
\nonumber
	\dot \phi_k &= \omega + K H_k\Big(R^{(0,\star)}(\phi) + K R^{(1,\star)}(\phi) + K^2 R^{(2,\star)}(\phi) + \mathcal O(K^3), \phi\Big)\\
	\nonumber
	&= \omega +K H_k(R_k^{(0,\star)}(\phi), \phi) + K^2 \nabla_R H_k(R^{(0,\star)}(\phi), \phi)\cdot R^{(1,\star)}(\phi) + \mathcal O(K^3)\\
	\label{eq:higher-order_expans}
	&= \omega +K H_k(1, \phi) + K^2 \nabla_R H_k(1, \phi)\cdot R^{(1,\star)}(\phi) + \mathcal O(K^3),
\end{align}
which is the base equation for phase reductions of any order. A phase reduction of first order can be obtained by truncating terms of order~$\mathcal O(K^2)$, a second order phase reduction is derived from~\eqref{eq:higher-order_expans} by ignoring all terms of order $\mathcal O(K^3)$, etc. In particular, the $(1,\infty)$-phase reduction is given by
\begin{align}
\nonumber
	\dot \phi_k &= \omega + KH_k(1,\phi)\\
	\label{eq:(1,inf)reduction}
	&= \omega + K\frac 1N \sum_{l=1}^N \left[ \frac{1+\delta g(\phi_l)}{1+\delta g(\phi_k)}\sin(\phi_l-\phi_k + \alpha) - \sin(\alpha)\right].
\end{align}
Up to now, this contains all orders of $\delta$, but by writing
\begin{align*}
	H_k(R,\phi) = H_k^{(-,0)}(R,\phi) + \delta H_k^{(-,1)}(R,\phi) + \delta^2 H_k^{(-,2)}(R,\phi) + \mathcal O(\delta^3),
\end{align*}
the $(1,\infty)$-phase reduction can also be written as
\begin{align*}
	\dot \phi_k \approx \sum_{n=0}^1\sum_{j = 0}^\infty K^n \delta^j P_k^{(n,j)}(\phi),
\end{align*}
where $P_k^{(0,0)}(\phi) \equiv \omega$, $P_k^{(0,j)}(\phi) \equiv 0$ for all $j\in \N$ and 
\begin{align*}
	P^{(1,j)}_k(\phi) = H_k^{(-,j)}(1,\phi), \quad j\in \N_0.
\end{align*}
For example, we find
\begin{align*}
	P_k^{(1,0)}(\phi) &= \frac 1N \sum_{l=1}^N \left[ \sin(\phi_l-\phi_k+\alpha)-\sin(\alpha)\right],\\
	P_k^{(1,1)}(\phi) &= \frac 1N \sum_{l=1}^N (g(\phi_l)-g(\phi_k))\sin(\phi_l-\phi_k+\alpha),\\
	P_k^{(1,2)}(\phi) &= \frac{1}{N} \sum_{l=1}^N g(\phi_k)(g(\phi_k)-g(\phi_l))\sin(\phi_l-\phi_k+\alpha).
\end{align*}
Consequently, the $(1,0)$-phase reduction is
\begin{align*}
	\dot \phi_k &= \sum_{n=0}^1 \sum_{j=0}^0 K^n \delta ^j P_k^{(n,j)}(\phi)\\
	&= \omega + \frac KN \sum_{l=1}^N \left[ \sin(\phi_l-\phi_k+\alpha) - \sin(\alpha) \right],
\end{align*}
which is the Kuramoto--Sakaguchi model~\cite{Sakaguchi1986} for identical oscillators.

\subsection{Higher-Order Phase Reductions}

Having explicitly stated the first-order phase reductions, we move on by considering the terms of order $K^2$ in~\eqref{eq:higher-order_expans}, thereby deriving a second-order phase reduction. As in Section~\ref{sec:prelim}, this requires knowledge of $R^{(1,\star)}(\phi)$. To get a formula describing it, we follow the lines of~\cite{Leon2019} and insert the expansion~\eqref{eq:R_expansion} into~\eqref{eq:main_system_R}. Using $R^{(0,\star)}\equiv 1$ and applying the chain rule, we find that the left-hand side of~\eqref{eq:main_system_R} turns into
\begin{align}
\nonumber
	\dot R_k &= \frac{\d}{\d t} R_k(t)\\
	\nonumber
	&= \frac{\d}{\d t}\Big( R_k^{(0,\star)}(\phi(t)) + K R_k^{(1,\star)}(\phi(t)) + \mathcal O(K^2)\Big)\\
	\nonumber
	&= K \nabla_\phi R_k^{(1,\star)}(\phi(t))\cdot \dot \phi(t) + \mathcal O(K^2)\\
	\nonumber
	&= K \nabla_\phi R_k^{(1,\star)}(\phi(t)) \cdot \Big( \omega \mathbbm 1 + K H(R, \phi)\Big) + \mathcal O(K^2)\\
	\label{eq:lhs}
	&= K\omega \nabla_\phi R_k^{(1,\star)}(\phi(t))\cdot \mathbbm 1 + \mathcal O(K^2),
\end{align}
whenever the dynamics is constrained to the limiting torus. Using, $\mathbbm{1} = (1,1,\dots, 1)^\top\in \R^N$ it follows that
\begin{align*}
	\omega \nabla \phi R_k^{(1,\star)}(\phi)\cdot \mathbbm{1} = \omega \sum_{l=1}^N \frac{\partial}{\partial \phi_l} R_k^{(1,\star)}(\phi).
\end{align*}
Similarly, the right-hand side turns into
\begin{align}
\nonumber
&F\Big(R_k^{(0)}(\phi(t)) + K R_k^{(1)}(\phi(t)) + \mathcal O(K^2)), \phi_k\Big) + KG_k\Big(R^{(0)}(\phi(t)) + K R^{(1)}(\phi(t)) + \mathcal O(K^2), \phi\Big)\\
\nonumber
&\quad = F(R_k^{(0)}(\phi),\phi_k) + K F_R(R_k^{(0)}(\phi), \phi_k) R_k^{(1)}(\phi) + K G_k(R^{(0)}(\phi), \phi) + \mathcal O(K^2)\\
\label{eq:rhs}
&\quad=F(R_k^{(0)}(\phi),\phi_k) + K\Big( F_R(R_k^{(0)}(\phi), \phi_k) R_k^{(1)}(\phi) + G_k(R^{(0)}(\phi), \phi)\Big) + \mathcal O(K^2),
\end{align}
where $F_R(R,\phi) = \frac{\partial}{\partial R} F(R,\phi)$. Now, equating~\eqref{eq:lhs} and~\eqref{eq:rhs} and collecting terms of order $\mathcal O(K^1)$ yields
\begin{align}\label{eq:R1PDE}
	F_R(1,\phi_k) R^{(1,\star)}_k(\phi) + G_k(1,\phi) = \omega \nabla_\phi R_k^{(1,\star)}(\phi)\cdot \mathbbm{1},
\end{align}
or equivalently, when using definitions of $F$ and $G$,
\begin{align*}
	&m ( 1+\delta g(\phi_k))^2 R^{(1,\star)}_k(\phi) + \frac{1}{N}\sum_{l=1}^N \Bigg[ \frac{1+\delta g(\phi_l)}{1+\delta g(\phi_k)}\cos(\phi_l-\phi_k+\alpha) - \cos(\alpha)\\
	& - \delta g'(\phi_k)\left(\frac{1+\delta g(\phi_l)}{(1+\delta g(\phi_k))^2}\sin(\phi_l-\phi_k+\alpha) - \frac{\sin(\alpha)}{1+\delta g(\phi_k)}\right) \Bigg] = \omega \nabla_\phi R_k^{(1,\star)}(\phi)\cdot \mathbbm 1,
\end{align*}
which is a linear first-order partial differential equation describing $R_k^{(1,\star)}(\phi)$.

At this point, we can no longer proceed as in Section~\ref{sec:prelim}, because we cannot set $\omega=0$, since our system is not rotationally invariant. Thus, we generalize the methods of~\cite{Leon2019} by solving the PDE~\eqref{eq:R1PDE}, as proposed in~\cite{Gengel2021}. Assuming an expansion
\begin{align}\label{eq:delta_expansion}
	R_k^{(1,\star)}(\phi) = R_k^{(1,0)}(\phi) + \delta R_k^{(1,1)}(\phi) + \delta^2 R_k^{(1,2)}(\phi) + \mathcal O(\delta^3),
\end{align}
we solve the PDE~\eqref{eq:R1PDE} order by order \cite{Evans2010, DeJager1996, Kevorkian1996}. When $\delta =0$, the PDE describing $R_k^{(1,0)}$ is 
\begin{align*}
	m R_k^{(1,0)}(\phi) + \frac{1}{N}\sum_{l=1}^N \left[ \cos(\phi_l-\phi_k+\alpha) - \cos(\alpha)\right] = \omega \nabla_\phi R_k^{(1,0)}(\phi)\cdot \mathbbm{1}.
\end{align*}
The solution to this PDE, which can be, for example, be found with the method of characteristics~\cite{Evans2010}, is given by
\begin{align}\label{eq:R10sol}
	R_k^{(1,0)}(\phi) = \frac{1}{Nm}\sum_{l=1}^N s_0(\phi_k,\phi_l), \qquad s_0(\phi_k, \phi_l) = -\cos(\phi_l-\phi_k + \alpha)  + \cos(\alpha).
\end{align}
On first-order in $\delta$, the resulting PDE is
\begin{align}
\label{eq:PDE_d1}
\begin{split}
	&m R^{(1,1)}_k(\phi)+ 2m g(\phi_k) R^{(1,0)}_k(\phi)- \omega \nabla_\phi R^{(1,1)}_k(\phi)\cdot \mathbbm{1}\\
	&\quad = -\frac{1}{N} \sum_{l=1}^N \Big[ (g(\phi_l)-g(\phi_k))\cos(\phi_l-\phi_k+\alpha) - g'(\phi_k)\left( \sin(\phi_l-\phi_k+\alpha)-\sin(\alpha)\right)\Big].
\end{split}
\end{align}
The solution of this PDE now depends on the specific choice of $g$. However, as one can infer from the structure of~\eqref{eq:PDE_d1}, its solutions are linear in $g$ in the sense that if $\hat R_k^{(1,1)}(\phi)$ is a solution to~\eqref{eq:PDE_d1} when $g = \hat g$ and $\tilde R_k^{(1,1)}(\phi)$ is one if $g = \tilde g$, then $\gamma \hat R_k^{(1,1)}(\phi)$ is a solution when $g = \gamma \hat g$ for all $\gamma \in \R$. Moreover, $\hat R_k^{(1,1)}(\phi) + \tilde R_k^{(1,1)}(\phi)$ is the solution to~\eqref{eq:PDE_d1} when $g = \hat g + \tilde g$.
If $g(\phi) = \sin(\phi)$, the solution of~\eqref{eq:PDE_d1} is given by
\begin{align}\label{eq:R11sol}
	R_k^{(1,1)}(\phi) = \frac{1}{2N(m^2+\omega^2)} \sum_{l=1}^N s_1(\phi_k, \phi_l),
\end{align}
where $s_1(\phi_k, \phi_l)$ is a trigonometric polynomial that is defined by
\begin{align}
\label{eq:s1}
\begin{split}
s_1(\phi_k, \phi_l) &= \omega \Big( 4\cos(\phi_l+\alpha) - \cos(\phi_k - 2\phi_l - \alpha) + 2\cos(2\phi_k-\phi_l-\alpha)\\
&\quad\qquad - 2\cos(\phi_k-\alpha) - 3 \cos(\phi_k+\alpha)\Big)\\
	&\quad +m \Big( 4\sin(\phi_l+\alpha) + \sin(\phi_k-2\phi_l-\alpha) + 2\sin(2\phi_k-\phi_l-\alpha)\\
	&\qquad\qquad - 2\sin(\phi_k-\alpha) - 3\sin(\phi_k+\alpha)\Big).
\end{split}
\end{align}
Equivalently, if $\alpha=0$, this can also be written as
\begin{align*}
	s_1(\phi_k, \phi_l) &= -2\Big(1-\cos(\phi_k-\phi_l)\Big)\Big(2\omega \cos(\phi_k)-\omega \cos(\phi_l) + 2m \sin(\phi_k)-m\sin(\phi_l)\Big)
\end{align*}
Finally, the PDE on order $\mathcal O(\delta^2)$ is
\begin{align*}
	&m R^{(1,2)}(\phi) + 2mg(\phi_k) R^{(1,1)}_k(\phi) +  m g(\phi_k)^2 R^{(0,1)}_k(\phi) - \omega \nabla_\phi R_k^{(1,2)}(\phi)\cdot \mathbbm 1\\
	& = - \frac 1N \sum_{l=1}^N \Bigg[ g'(\phi_k)\Big( (g(\phi_k)-g(\phi_l))\sin(\phi_l-\phi_k+\alpha) + g(\phi_k)(\sin(\phi_l-\phi_k+\alpha)-\sin(\alpha))\Big)\\
	& \qquad \qquad \qquad+ g(\phi_k)(g(\phi_k)-g(\phi_l)) \cos(\phi_l-\phi_k+\alpha) \Bigg].
\end{align*}
This PDE, however, is not linear in $g$. In particular, if $\hat R^{(1,2)}_k(\phi)$ is a solution for $g=\hat g$ then $\gamma^2 \hat R^{(1,2)}_k(\phi)$ solves the PDE whenever $g = \gamma \hat g$, for $\gamma\in \R$. If $g(\phi) = \sin(\phi)$ a solution is of the form
\begin{align*}
	R_k^{(1,2)}(\phi) &= \frac{1}{-4 m N \left(m^4+5 m^2 \omega ^2+4 \omega ^4\right)} \sum_{l=1}^N s_2(\phi_k, \phi_l),
\end{align*}
where $s_2(\phi_k, \phi_l)$ is a trigonometric polynomial of the same form as $s_1$ but with more summands\footnote{The full expression of $s_2(\phi_k, \phi_l)$ can be generated with the \textsc{Mathematica} code accompanying the paper, see~\cite{Bohle2023GitHub}.}.

To determine the second order interactions in~\eqref{eq:higher-order_expans}, we also need to expand $\nabla_R H_k(R^{(0)}(\phi), \phi)$ in terms of $\delta$. Denoting $H_k(R,\phi) = H_k^{(-,0)}(R,\phi) + \delta H_k^{(-,1)}(R,\phi) + \delta^2 H_k^{(-,2)}(R,\phi) + \mathcal O(\delta^3)$, we find
\begin{align*}
	\nabla_R H_k^{(-,0)}(1,\phi) &=\frac 1N  \begin{pmatrix} \sin(\phi_1-\phi_k + \alpha) \\ \vdots \\ \sin(\phi_N-\phi_k + \alpha)
	\end{pmatrix}
	 - \frac 1N e_k \sum_{l=1}^N \sin(\phi_l-\phi_k+\alpha),\\
	 \nabla_R H_k^{(-,1)}(1, \phi) &= \frac{1}{N} \begin{pmatrix}
	 (g(\phi_1)-g(\phi_k))\sin(\phi_1-\phi_k+\alpha)\\ \vdots \\ (g(\phi_N)-g(\phi_k))\sin(\phi_N-\phi_k+\alpha)
	 \end{pmatrix}\\
	 &\qquad -\frac{1}{N}e_k \sum_{l=1}^N (g(\phi_l)-g(\phi_k)) \sin(\phi_l-\phi_k+\alpha),\\
	 \nabla_R H^{(-,2)}_k(1, \phi) &= \frac{1}{N}\begin{pmatrix}
	 g(\phi_k)(g(\phi_k)-g(\phi_1))\sin(\phi_1-\phi_k+\alpha)\\ \vdots \\ g(\phi_k)(g(\phi_k)-g(\phi_N))\sin(\phi_N-\phi_k+\alpha)
	 \end{pmatrix}\\
	 &\quad -\frac{1}{N}e_k \sum_{l=1}^N  g(\phi_k) (g(\phi_k)-g(\phi_l)) \sin(\phi_l-\phi_k+\alpha),
\end{align*}
where $e_k$ is the $k$th unit vector in $\R^N$.
Finally, we can put everything together and calculate the second order terms in~\eqref{eq:higher-order_expans}:
\begin{align*}
\nabla_R H_k(1, \phi) \cdot R_k^{(1,\star)}(\phi) = P^{(2,0)}_k(\phi) + \delta P^{(2,1)}_k(\phi) + \delta^2 P^{(2,2)}_k(\phi) + \mathcal O(\delta^3),
\end{align*}
with
\begin{subequations}
\label{eq:P2}
\begin{align}
	P^{(2,0)}_k(\phi) &= \nabla_R H_k^{(-,0)}(1,\phi) \cdot R^{(1,0)}(\phi),\\
\label{eq:P21}
	P^{(2,1)}_k(\phi) &= \nabla_R H_k^{(-,0)}(1,\phi) \cdot R^{(1,1)}(\phi) + \nabla_R H_k^{(-,1)}(1,\phi)\cdot R^{(1,0)}(\phi),\\
	\begin{split}
	P^{(2,2)}_k(\phi) &= \nabla_R H_k^{(-,0)}(1,\phi)\cdot R^{(1,2)}(\phi) + \nabla_RH_k^{(-,1)}(1,\phi)\cdot R^{(1,1)}(\phi)\\
	&\qquad + \nabla_R H_k^{(-,2)}(1,\phi)\cdot R^{(1,0)}(\phi).
	\end{split}
\end{align}
\end{subequations}
Evaluating these expressions yields, for example, 
\begin{align*}
	P^{(2,0)}_k(\phi) = \frac{1}{2N^2m}\sum_{l=1}^N \sum_{i=1}^N &\Big( \sin(\phi_i+\phi_k-2\phi_l)- \sin(\phi_i-\phi_k+2\alpha) + \sin(\phi_i - 2 \phi_k + \phi_l + 2 \alpha) \Big)
\end{align*}
for the second-order terms in $K$ when $\delta = 0$. We emphasize that it is here clearly visible that three phases interact with each other, which is different from the terms that appear in a $(1,0)$ phase reduction.

In conclusion, the $(2,2)$-phase reduction is given by
\begin{align*}
	\dot \phi_k &= \omega + K\Big( P_k^{(1,0)}(\phi) + \delta P_k^{(1,1)}(\phi) + \delta ^2 P_k^{(1,2)}(\phi)\Big)\\
	&+ K^2 \Big( P_k^{(2,0)}(\phi) + \delta P_k^{(2,1)}(\phi) + \delta ^2 P_k^{(2,2)}(\phi)\Big),
\end{align*}
with $P_k^{(1,0)}(\phi), P_k^{(1,1)}(\phi)$ and $P_k^{(1,2)}(\phi)$ are as in Section~\ref{sec:first_order} and $P_k^{(2,0)}(\phi), P_k^{(2,1)}(\phi)$ and $P_k^{(2,2)}(\phi)$ are defined in~\eqref{eq:P2}.

\subsection{Comparison of Phase Reductions with and without Symmetry}

As we have highlighted in this section, the full system~\eqref{eq:main_system} has a rotational $\S$~symmetry for $\delta=0$ that breaks for a generic perturbation to the limit cycle and $\delta \neq 0$.

We now consider the reduced phase equations in detail focusing on the corresponding change in symmetry properties one would expect\footnote{Note that the phase reduction is in terms of the original oscillator phases to analyze the symmetry properties explicitly. In particular, we do not consider an additional near-identity transformation of the phases to put the phase equations in normal form as in~\cite{VonderGracht2023a} or an additional averaging approximation that can lead to more symmetries in the reduction than in the full nonlinear system; cf.~\cite{Crawford1991}.}.
The $(2,0)$-phase reduction, i.e., the second order phase reduction when $\delta = 0$, is given by
\begin{align*}
	\dot \phi_k &= \omega + K P_k^{(1,0)}(\phi) + K^2 P_k^{(2,0)}(\phi)\\
	&= \omega + K \frac 1N \sum_{l=1}^N[ \sin(\phi_l - \phi_k + \alpha) - \sin(\alpha)]\\
	&\qquad + K^2 \frac{1}{2N^2m}\sum_{l=1}^N \sum_{i=1}^N \Big( \sin(\phi_i+\phi_k-2\phi_l)- \sin(\phi_i-\phi_k+2\alpha)\\
	&\qquad \qquad \qquad\qquad \qquad \quad+ \sin(\phi_i - 2 \phi_k + \phi_l + 2 \alpha) \Big).
\end{align*}
As one can see, the right-hand side of this equation only depends on phase differences. Therefore, its value remains invariant when shifting all oscillators by a common phase. Consequently, the $(2,0)$-phase reduction inherits the $\S$ symmetry of the full system. 
Writing $Q e^{i\Theta} = \frac 1N \sum_{j=1}^N e^{2 i \phi_j}$ and $R e^{i\Psi} = \frac 1N \sum_{j=1}^N e^{i\phi_j}$ as in~\cite{Leon2019}, we can compare the results in~\cite{Leon2019} with our $(2,0)$-phase reduction. 
With this notation the $(2,0)$-phase reads
\begin{align*}
	\dot \phi_k &= \omega + K R \sin(\Psi - \phi_k + \alpha) - K \sin(\alpha)\\
	&\qquad + K^2\frac{1}{2m}\left(RQ\sin(\Psi + \phi_k - \Theta) - R\sin(\Psi - \phi_k + 2\alpha)+R^2\sin(2\Psi - 2\phi_k+2\alpha)\right),
\end{align*}
which agrees with the result from~\cite[Equation (15)]{Leon2019}, if one chooses $K=\epsilon|1+ic_1|$ and $m = -2$. 
This shows that even though the nonlinearity in the radial direction in~\eqref{eq:single_oscillator_r} is different from the nonlinearity of the Stuart--Landau oscillator considered in~\cite{Leon2019}, the phase equations up to second order in~$K$ agree as expected.
The plausibility of $K=\epsilon | 1+ic_1|$ is immediate, when comparing the coupling strength in~\eqref{eq:LP_CGL} with those in~\eqref{eq:coupling}. Moreover, $m=-2$ can be explained as follows: When $\delta = 0$, $m$ is the rate of attraction towards the limit cycle of a single oscillator~\eqref{eq:single_oscillator}. In particular, this rate can be obtained by linearizing~\eqref{eq:single_oscillator_r} with respect to $r$ and evaluating at the limit cycle $r=1$. Doing the same for the oscillator~\eqref{eq:LP_single}, yields that the rate of attraction to this limit cycle is $-2$.

Now, let us consider phase reductions when $\delta\neq 0$. When $g(\phi) = \sin(\phi)$, the $(2,1)$-phase reduction is given by the system
\begin{align*}
	\dot \phi_k &= \omega + K P_k^{(1,0)}(\phi)+ K\delta P_k^{(1,1)}(\phi)+ K^2 P_k^{(2,0)}(\phi)+K^2 \delta P_k^{(2,1)}(\phi)\\
	&=  \omega + K \frac 1N \sum_{l=1}^N[ \sin(\phi_l - \phi_k + \alpha) - \sin(\alpha)] + K\delta \frac 1N \sum_{l=1}^N (\sin(\phi_l)-\sin(\phi_k))\sin(\phi_l-\phi_k+\alpha) \\
	&\qquad + K^2 \frac{1}{N^2m}\sum_{l=1}^N \sum_{i=1}^N \sin(\phi_l-\phi_k+\alpha)( -\cos(\phi_i-\phi_l+\alpha)+\cos(\phi_i-\phi_k+\alpha) )\\
	&\qquad + K^2\delta \frac{1}{2N^2(m^2+\omega^2)}\sum_{l=1}^N\sum_{i=1}^N \sin(\phi_l-\phi_k + \alpha) (s_1(\phi_l, \phi_i)-s_1(\phi_k,\phi_i)) \\
	&\qquad + K^2\delta \frac{1}{N^2m} \sum_{l=1}^N\sum_{i=1}^N\Big( (\sin(\phi_l)-\sin(\phi_k))\sin(\phi_l-\phi_k+\alpha)\\
	&\qquad \qquad \qquad \qquad \qquad \quad \cdot (-\cos(\phi_i - \phi_l+\alpha) + \cos(\phi_i-\phi_k+\alpha))\Big),
\end{align*}
which does not have an~$\S$~symmetry. In particular, the terms of order $\delta$ are not invariant when one shifts all oscillators by a common phase. Since, higher-order phase reductions also consist of these terms, any phase reduction of higher-order is not $\S$ symmetric. To conclude, phase reductions of the full system~\eqref{eq:main_system} possess an $\S$ symmetry if and only if the full system itself possesses this symmetry.

As a remark, when $\alpha = 0$, the $(2,1)$-phase reduction can also be written as
\begin{align*}
	\dot \phi_k &= \omega + K\frac 1N \sum_{l=1}^N \sin(\phi_l-\phi_k)\Big[ 1 + \delta (\sin(\phi_l)-\sin(\phi_k))\Big]\\
	&\qquad + K^2 \frac{1}{N^2m} \sum_{l=1}^N\sum_{i=1}^N \Bigg\{ \sin(\phi_l-\phi_k)(1-\cos(\phi_i-\phi_l)) \left[ 1 + \delta\big( u(\phi_l, \phi_i)+ \sin(\phi_l)-\sin(\phi_k)\big) \right] \\
	&\qquad\qquad\qquad\qquad\qquad - \sin(\phi_l-\phi_k)(1-\cos(\phi_i-\phi_k)) \left[1 + \delta\big( u(\phi_k, \phi_i)+ \sin(\phi_l)-\sin(\phi_k) \big) \right] \Bigg\},
\end{align*}
where
\begin{align*}
	u(\phi_k, \phi_i) = \frac{-m}{m^2+\omega^2} (2\omega \cos(\phi_k) - \omega \cos(\phi_i) + 2m \sin(\phi_k)-m\sin(\phi_i)).
\end{align*}
This formula shows that effect of $\delta$ on the $(2,1)$-phase reduction can also be seen as a perturbation of the $(2,0)$-phase reduction.

\section{Dynamics}\label{sec:dynamics}

In this section we consider two different orbits, i.e., the synchronized orbit and the splay orbit, and compare their stability in a few different systems, including the full system~\eqref{eq:main_system} and various phase reductions on different order.

\subsection{Synchronized State}\label{sec:dynamics_sync}

First, we consider the synchronized state in the full system. A state in the phase space is called synchronized if all the oscillators are at the same position. In the full system~\eqref{eq:main_system} this state is defined by $\{A_1=\dots = A_N\}$ or in polar coordinates 
\begin{align}\label{eq:sync_state}
	\{R_1 = \dots = R_N\} \quad \text{and} \quad \{\phi_1=\dots = \phi_N\}.
\end{align}
Consequently, if a state is synchronized, it is uniquely given by its amplitude $R^\star := R_i$ and its phase $\phi^\star := \phi_i$ for any $i=1,\dots,N$. Due to the $S_N$ symmetry of~\eqref{eq:main_system} the set of synchronized states~\eqref{eq:sync_state} is dynamically invariant. Thus, we can insert the ansatz~\eqref{eq:sync_state} into the system~\eqref{eq:main_system} to obtain ODEs for $R^\star$ and $\phi^\star$ as
\begin{align*}
	\dot R^\star &= m (R^\star)^2 (R^\star - 1)(1+\delta g(\phi^\star))^2,\\
	\dot \phi^\star &= \omega.
\end{align*}
Based on these equations one can see that always $R^\star = 1$ on the invariant torus and that the rate of attraction to $R^\star = 1$ is given by $m (1+\delta g(\phi^\star))^2$. Since $\abs{\delta}$ is small, this rate is mostly governed by $m<0$. Moreover, the phase $\phi^\star$ evolves with constant speed $\omega$. Consequently, the synchronized orbit on the invariant torus is given by $\gamma^\mathrm f(t)=(\mathbbm 1, (\omega_0 + \omega t)\mathbbm 1)$, which is periodic with period $T=2\pi/\omega$.

Now, we look at synchronized states in phase-reduced systems. In a phase-reduced system, there are no amplitudes and thus, a synchronized state is present when the single condition $\phi_1 = \dots = \phi_N$ is fulfilled. Here, a synchronized state is only determined by its phase $\phi^\star := \phi_i$ for any $i=1,\dots,N$. Similarly to the full system, phase-reduced systems retain the $S_N$ symmetry and therefore the set of synchronized states is dynamically invariant. Inserting $\phi\equiv \phi^\star$ into any of the phase reductions derived in Section~\ref{sec:meat} yields
\begin{align*}
	\dot \phi^\star = \omega.
\end{align*}
Therefore, the synchronized orbit in phase-reduced systems is $\gamma^\mathrm{pr}(t) = (\omega_0 + \omega t)\mathbbm 1$ with period $T=2\pi /\omega$.

Having established representations for the synchronized orbits, we now investigate their stability. Usually, when one wants to check the stability of a synchronized orbit, one changes into a co-rotating system, in which each synchronized state is an equilibrium. Then, one linearizes the vector field around this equilibrium and calculates the eigenvalues of this linearization. If they are all negative, apart from a single $0$ eigenvalue that corresponds to perturbations along the continuum of synchronized states, the synchronized orbit is linearly stable and linearly (neutrally) unstable otherwise. However, when $\delta\neq 0$, we cannot change into a co-rotating coordinate system, since the full system as well as phase-reduced systems are not $\S$ symmetric. In particular, the rate of attraction to the limit cycle depends on the position on the limit cycle. Consequently, one needs to take averages over all rates of attraction of one period of the limit cycle. These averages are called Floquet exponents. The concept of Floquet exponents and Poincar\'e return maps is often helpful when analyzing the stability of periodic orbits \cite{Chicone2006, Teschl2012}. To understand it, let us consider the general differential equation
\begin{align*}
	\dot x = \mathcal H(x), \qquad \mathcal H\colon \mathcal X\to \textnormal{T}\mathcal X,
\end{align*}
where $\mathcal H$ is a smooth vector field on the phase space $\mathcal{X}$ and $\textnormal{T}\mathcal X$ is the tangent bundle. Suppose there is a periodic orbit $\gamma\colon [0,T]\to \mathcal X$ with period~$T$. Later we want to apply this to the full system, in which $x=(R,\phi)^\top\in \R_{\ge 0}^N\times \S^N$, and phase-reduced systems with $x=\phi\in \S^N$. To analyze the stability of the orbit $\gamma$, one assumes a perturbation of the starting point of the periodic orbit $x(0) = \gamma(0) + \epsilon \eta(0)$ and continues this perturbation along the periodic orbit such that $x(t) = \gamma(t) + \epsilon\eta(t)$. One the one hand, if all possible perturbations $\eta(0)$ decay after one period $T$, we expect the periodic orbit $\gamma$ to be stable. On the other hand, if some perturbations $\eta(0)$ grow in amplitude, the orbit $\gamma$ is unstable. Therefore, we solve for $\eta(t)$. However, since that is difficult to do in general, we first linearize in~$\epsilon$ to obtain on first order
\begin{align}\label{eq:Floquet_orbit}
\dot \eta(t) = A(t)\eta(t),\quad A(t) = D\mathcal H(\gamma(t)),
\end{align}
which is a system of linear ODEs with a time-dependent coefficient matrix $A(t)$. Now, let $\Phi(t)\in \R^{N\times N}$ be a fundamental solution of this ODE such that $\eta(t) = \Phi(t)\eta(0)$. To determine the linear stability of the periodic orbit $\gamma(t)$ one propagates all possible perturbations $\eta(0)$ over one period $T$ of the orbit and then looks at the eigenvalues of the map $\Phi(T)$. Of course, this map has an eigenvalue $1$, that corresponds to the eigenvector that represents a perturbation along the periodic orbit. All other eigenvalues $\lambda_k$, $k=1,\dots,N-1$ are the multipliers of a Poincar\'e return map. They are related to the Floquet exponents $q_k$ of the orbit as $\lambda_k = e^{T q_k}$, for $k=1,\dots, N-1$. If the largest absolute value of all Poincar\'e return map multipliers (PRMMs) is less than $1$, the periodic orbit is linearly stable. If one PRMM has an absolute value greater than $1$, the orbit is linearly unstable.

Next, we apply this concept to the synchronized orbit in phase-reduced systems. We start by calculating the PRMMs for the $(1,\infty)$-phase-reduced system~\eqref{eq:(1,inf)reduction}. When putting this system into the framework of~\eqref{eq:Floquet_orbit}, we see that the matrix $A(t)$ is 
\begin{align*}
	A(t) = f(\gamma(t))\frac{K}{N}( \mathbbm 1 - N \mathbb I)
\end{align*}
with
\begin{align*}
	f(\gamma) = \frac{1}{1+\delta g(\gamma)}\Big( \delta g'(\gamma) \sin(\alpha) + (1+\delta g (\gamma))\cos(\alpha)\Big),
\end{align*}
where $\mathbbm 1$ is a $N\times N$ matrix where all entries are ones and $\mathbb I$ is the $N\times N$ dimensional identity matrix. Due to the $S_N$ symmetry of the synchronized state in the phase-reduced systems, the matrix $A(t)$ has the special property that it is just a multiple of $\mathbbm 1-N\mathbb I$. Since matrices of this form commute with each other, a fundamental solution $\Phi(t)$ of~\eqref{eq:Floquet_orbit} can explicitly be calculated using the matrix exponential:
\begin{align*}
	\Phi(t) &= \exp\left(  \int_0^t A(\hat t) \mathrm d \hat t  \right)\\
	&=\exp\left( \int_0^t f(\gamma(\hat t)) \mathrm d \hat t \left( \frac KN (\mathbbm 1 - N \mathbb I)\right) \right).
\end{align*}
Integrating $f(\gamma(t))$ over one period of the orbit $\gamma(t) = \omega_0 + \omega t$ yields
\begin{align*}
	\int_0^T f(\gamma(\hat t)) \ \mathrm d \hat t &= \left[\sin(\alpha)\ln(1+\delta g(\gamma(\hat t))) + \hat t\cos(\alpha)\right]_{\hat t = 0}^{\hat t = T}\\
	&=T\cos(\alpha) = \frac{2\pi}{\omega} \cos(\alpha).
\end{align*}
Combining this with the fact that the eigenvalues of $\exp( \frac aN (\mathbbm 1 -N\mathbb I) )$ are $e^{-a}$ with multiplicity $N-1$ and $1$ with multiplicity $1$, we infer that the critical PRMM is
\begin{align*}
	\lambda^\mathrm{crit} = \exp \left( \frac{-2\pi K}{\omega}\cos(\alpha) \right).
\end{align*}
Interestingly, this is independent of $\delta$ even though the system $(1,\infty)$-phase reduction includes all orders in $\delta$. Consequently, the stability of a synchronized orbit is unaffected by $\delta$ in a phase reduction without higher-order interactions. A similar calculation yields that the critical PRMM of the synchronized orbit in a $(2,0)$-phase reduction is
\begin{align}\label{eq:floquet_crit_(2,0)}
	\lambda^\mathrm{crit} = \exp \left( \frac{-2\pi K}{m\omega}(m\cos(\alpha) - K\sin(\alpha)^2)\right).
\end{align}
The critical PRMM in a $(2,1)$-phase reduction agrees with~\eqref{eq:floquet_crit_(2,0)}, which can be shown using the formulas~\eqref{eq:R11sol} and~\eqref{eq:P21} if $g(\phi) = \sin(\phi)$. 
The independence of $\lambda^\mathrm{crit}$ on~$\delta$ in a $(2,1)$-phase reduction can be explained by a symmetry mapping~$\delta$ to~$-\delta$. Assuming $g(\phi)=\sin(\phi)$, the symmetry would alter the system such that the limit cycle is then parameterized by $r=1-\delta g(\phi)$. Using the $\S$-symmetry of the system one can apply a phase shift $\phi \mapsto \phi + \pi$ to obtain the original system where the limit cycle is given by $r=1-\delta g(\phi+\pi)=1+\delta g(\phi)$. Consequently, the $\S$-symmetry causes a~$\Z_2$ symmetry mapping~$\delta$ to~$-\delta$. Contributions of~$\delta$ to $\lambda^\mathrm{crit}$ can thus only be via even powers of $\delta$, which is why it~$\delta$ does not appear in the critical PRMM of a $(2,1)$-phase reduction. When $g$ is a general function, $\lambda^\mathrm{crit}$ of a $(2,1)$-phase reduction still does not depend on~$\delta$ and the formulas to show that are found in the Appendix. The stability of the synchronized orbit in a $(2,1)$-phase reduction thus agrees with its stability in a $(2,0)$-phase reduction, i.e., one for which the limit cycle is circular. Thus, when concerned with the stability of the synchronized orbit, deformations of the limit cycle can be ignored to first order.
Finally, when $g(\phi) = \sin(\phi)$, the critical PRMM in a $(2,2)$-phase reduction it is given by
\begin{align*}
	\lambda^\mathrm{crit} & = \exp\Bigg( \frac{-2\pi K}{m\omega (m^2+\omega^2)}\Big( m^3 \cos(\alpha) + m\omega^2 \cos(\alpha) - Km^2\sin(\alpha)^2\\
	&\qquad \qquad \qquad \qquad \qquad - 2Km^2\delta^2\sin(\alpha)^2 - K\omega^2\sin(\alpha)^2 \Big)\Bigg),
\end{align*}
which finally shows the dependence on $\delta$.

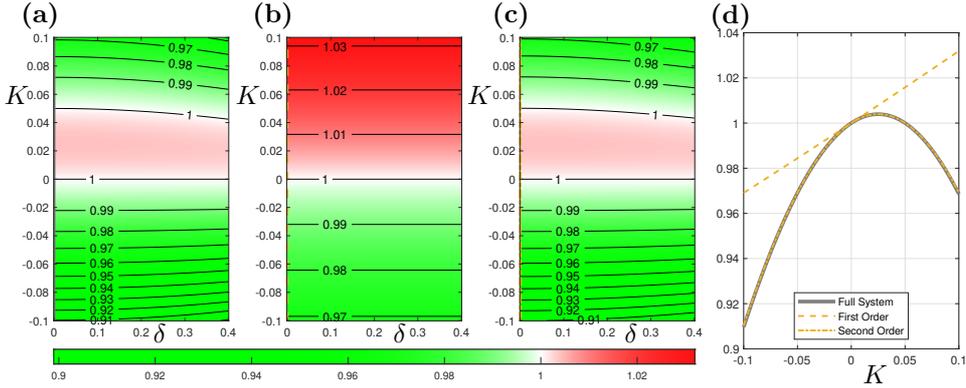
\begin{figure}[ht]
\begin{overpic}[width = 0.74\textwidth, grid = false]{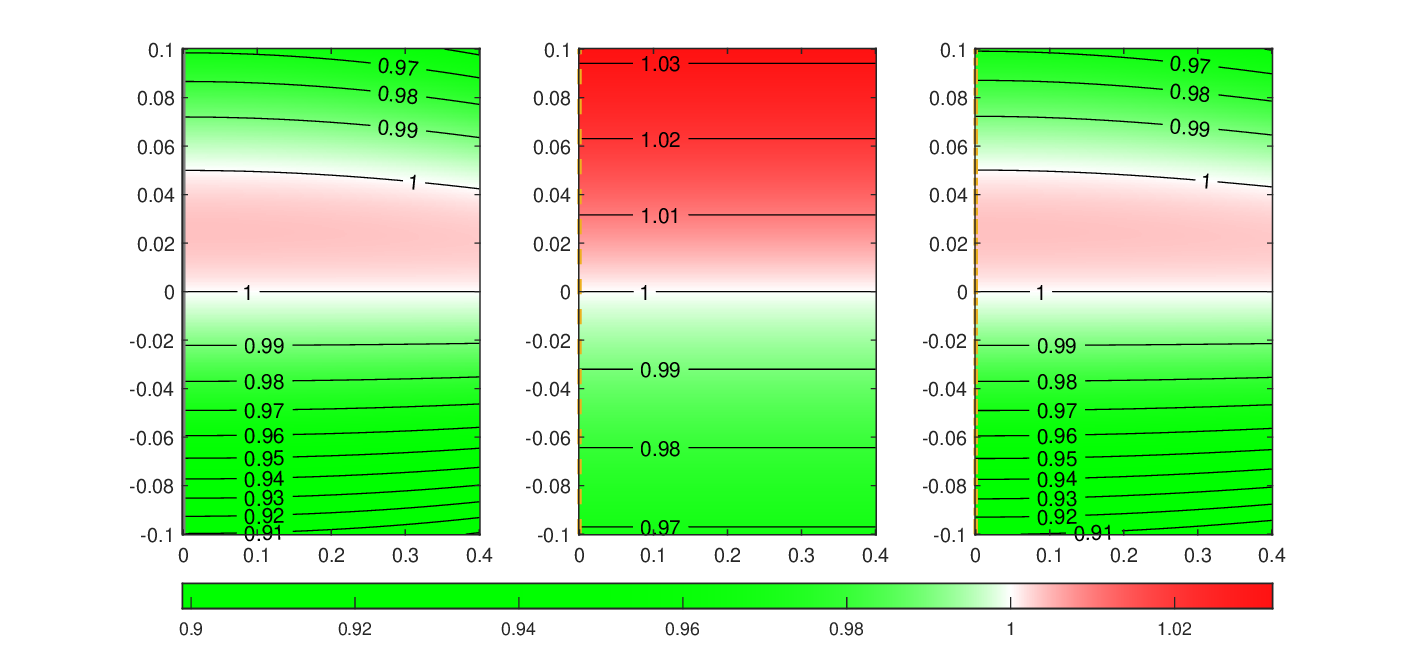}
\put(9,45){\textbf{(a)}}
\put(37,45){\textbf{(b)}}
\put(66,45){\textbf{(c)}}
\put(25,6){$\delta$}
\put(54,6){$\delta$}
\put(82,6){$\delta$}
\put(7,35){$K$}
\put(35,35){$K$}
\put(63,35){$K$}
\end{overpic}\hspace{-1cm}
\begin{overpic}[width = .25\textwidth, grid=false]{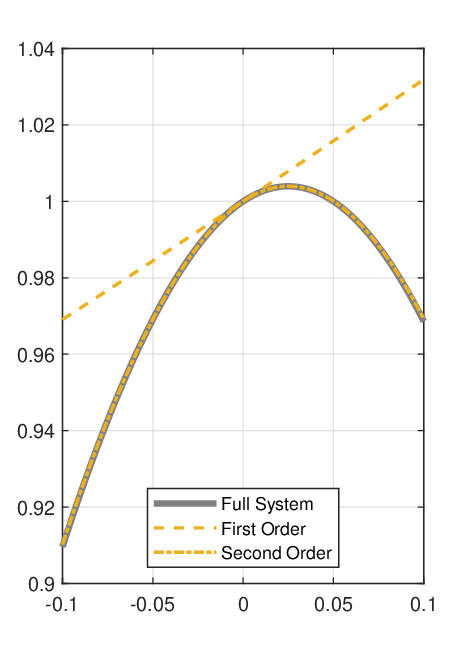}
\put(1,95){\textbf{(d)}}
\put(40,2){$K$}
\end{overpic}

\caption{Comparison of the critical multipliers of the Poincar\'e return map of the synchronized orbit in the full system and phase-reduced systems. (a) shows the critical PRMM of the synchronized orbit in the full system in dependence of $\delta$ and the coupling strength $K$. Part (b) depicts the critical PRMM of the synchronized orbit in a $(1,\infty)$-phase reduction and part (c) illustrates the critical PRMM of the same orbit in a $(2,2)$-phase reduction. 
Finally (d) depicts the critical PRMM in the full system, the $(1,\infty)$-phase reduction and the $(2,2)$-phase reduction when $\delta=0$. Parameter values: $\omega = 1, m=-1, \alpha = \pi/2+1/20, g(\phi) = \sin(\phi)$.}
\label{fig:floquet_sync}
\end{figure}

Having derived stability conditions for synchronized orbits in phase-reduced system, we now analyze the stability of the synchronized orbit in the full system. Unfortunately, when applying the concept~\eqref{eq:Floquet_orbit} to the full system, the matrices $A(t)$ and $A(s)$ do not commute with each other. Therefore, it is not possible to use the matrix exponential to analytically compute PRMMs or Floquet exponents, but we need to use numerical methods to determine them. 
Yet, in the special case $\delta = 0$, the stability analysis of these periodic orbits simplifies quite significantly. In fact, in this case, the full system has an $\S$~symmetry, which acts by shifting all oscillators by a constant phase. Then, one can also change to a co-rotating coordinate frame, in which $\omega = 0$. In these new coordinates all synchronized states are then equilibria. The spectrum of the right-hand side of the system at the synchronized state then contains information about the stability. There will be one $0$ eigenvalue, since there is a one-dimensional continuum of synchronized states. If all the other eigenvalues have negative real part, the synchronized state as an equilibrium in the co-rotating frame is linearly stable and thus the synchronized orbit in the original system inherits this stability. Conversely, if one eigenvalue has positive real part, the synchronized orbit in the original system is unstable. Conducting this analysis for the full system~\eqref{eq:main_system} yields that the linearization of the right-hand side is given by a matrix
\begin{align}\label{eq:full_system_Drhs}
	\begin{pmatrix}
		m \mathbb I + \frac KN \cos(\alpha) (\mathbbm 1 - N \mathbb I)& -\frac{K}{N}\sin(\alpha) (\mathbbm 1 - N\mathbb I)\\
		\frac{K}{N}\sin(\alpha) (\mathbbm 1 - N\mathbb I) & \frac{K}{N}\cos(\alpha)(\mathbbm 1 - N\mathbb I)
	\end{pmatrix}.
\end{align}
The eigenvalues of this matrix can explicitly be calculated and are given by
\begin{align*}
	q_1 &= 0,\\
	q_{2,\dots,N} &= \frac 12 \left( m - 2K\cos(\alpha) + \sqrt{-2K^2+m^2+2K^2\cos(2\alpha)} \right),\\
	q_{N+1} &= m,\\
	q_{N+2,\dots,2N} &= \frac 12 \left( m - 2K\cos(\alpha) - \sqrt{-2K^2+m^2+2K^2\cos(2\alpha)} \right),
\end{align*}
where we denote them by $q_k$ for $k=1,\dots, 2N$ because they describe the instantaneous rate of attraction to the periodic orbit and thus relate to the Floquet exponents. In fact, since this instantaneous rate of attraction is constant over the whole orbit, these eigenvalues agree with the Floquet exponents.
While the first~$N$ eigenvalues correspond to perturbations of the phases~$\phi$, the last $N$ eigenvalues originate from perturbations in the radial directions. To compare this model with phase-reduced models, we assume that $-m$ is big enough such that the last $N$ eigenvalues can be neglected and the critical Floquet exponent $q^\textrm{crit}$ is given by $q_2,\dots,q_N$. Of course there is also the zero eigenvalue~$q_1$. However, that does not contribute to the stability as it corresponds to a perturbation along the continuum of synchronized states. Given the critical Floquet exponent, one can then obtain the critical PRMM by simply calculating $\lambda^\textrm{crit} = \exp( q^\textrm{crit} T) = \exp( \frac{2\pi}{\omega} q^\textrm{crit})$.
When $\delta \neq 0$ in the full system, PRMMs can only numerically be calculated, as shown in Figure~\ref{fig:floquet_sync}(a). Figures~\ref{fig:floquet_sync}(b-d) compare PRMMs from the full system with PRMMs from phase-reduced systems.

\subsection{Splay State}

While all oscillators gather at one point on the circle if they are synchronized, one can say that the splay state is the opposite of that. A splay state is given when the oscillators phases are equidistantly distributed on the circle. More specifically, a state $\phi\in\S^N$ is a splay state if there is a permutation $\sigma\colon [N]\to [N]$ such that
\begin{align*}
	\phi_{\sigma(k+1)} = \phi_{\sigma(k)} + \frac{2\pi}{N},
\end{align*}
for $k=1,\dots,N-1$. By relabeling the nodes, we might also assume that $\sigma$ is the identity map. Thus the set of splay states is given by
\begin{align}\label{eq:splay}
	D := \{ \phi\in \S^N: \phi_{k+1} = \phi_k + \frac{2\pi}{N} \text{ for } k = 1,\dots,N-1\}.
\end{align}
A splay state in this set can be characterized by just the first phase $\phi_1$. This set has a~$Z_N := \Z/(N\Z)$ symmetry group that acts on a state by shifting the indices, which have to be understood modulo~$N$, of each oscillator by a constant integer. Now, let us consider the set of splay states in phase-reduced systems with $\delta = 0$. Since the right-hand sides of $(1,0)$- and $(2,0)$-phase-reduced systems are equivariant with respect to this group action, the set of splay states is dynamically invariant. In particular, when inserting a splay state into the right-hand side of a $(1,0)$-phase reduction and a $(2,0)$-phase reduction, it follows that $\dot \phi_k = \omega - K\sin(\alpha)=:\hat \omega$. Therefore, if $\hat \omega\neq 0$, there exists a periodic orbit $\gamma(t)\in \S^N$ with $\gamma_k(t) := \hat\omega t + \frac{2\pi k}{N}$, that has period $T = \frac{2\pi}{\hat \omega} = \frac{2\pi}{\omega-K\sin(\alpha)}$. We refer to this orbit as the splay orbit.

Next, we consider splay states in the full system~\eqref{eq:main_system} and still assume $\delta = 0$. Inserting the ansatz~\eqref{eq:splay} into the full system~\eqref{eq:main_system} yields that the amplitudes on the invariant torus are given by 
\begin{align}\label{eq:splay_amplitude}
	R_k \equiv\frac 12 \left(1+\sqrt{1+\frac{4K\cos(\alpha)}{m}}\right) =: R^\star.
\end{align}
Therefore, we call a state $(R,\phi)\in \R_{\ge 0}^N\times \S^N$ a splay state if $\phi\in D$ and the amplitudes satisfy~\eqref{eq:splay_amplitude}. Then, the splay state in the full system has the same symmetry group $Z_N$ as splay states in phase-reduced systems. Moreover, since the right-hand side of the full system~\eqref{eq:main_system} is again equivariant with respect to this symmetry group, the splay state is dynamically invariant. Furthermore, the angular frequency of the phases is $\dot \phi_k = \omega - K\sin(\alpha) =\hat \omega$ as for phase-reduced systems. Therefore, there there exists a periodic orbit $\gamma(t) = (R(t), \phi(t))^\top$ with $R_k(t) = R^\star$ and $\phi_k(t) = \hat \omega t + \frac{2\pi k }{N}$, that has the same period as the one in phase-reduced systems.

When analyzing the stability of these splay orbits in both phase-reduced and the full system, it is important to note that splay orbits are just one single periodic orbit in a whole continuum of periodic orbits. In particular, in phase-reduced systems, all incoherent states, that are characterized by
\begin{align}\label{eq:incoherent}
	Z := \frac 1N \sum_{k=1}^N e^{i\phi_k} = 0,
\end{align}
rotate around the circle with constant frequency~$\hat \omega$ and are a union of manifolds~\cite{Ashwin2016x}.
In the full system states $(R,\phi)\in \R_{\ge 0}^N\times \S^N$ with $R_k= R^\star$ for $k=1,\dots,N$ and $Z=0$ rotate around the circle with the same constant frequency. Therefore, there exist further periodic orbits, which we refer to as incoherent orbits. Since a splay state is a special incoherent state, but in general the set of incoherent states is larger than the set of splay states, the splay orbit is only one orbit in a continuum of incoherent orbits. Consequently, when analyzing the stability of the splay orbits with PRMMs or Floquet exponents, there is always one neutral multiplier or exponent, respectively. The only exception is present when the set of incoherent states coincides with the set of splay states, i.e., when $N=3$. To overcome the problem of neutral stability, we restrict ourselves to $N=3$.

If $\delta = 0$, we can change into a rotating frame coordinate system, in which splay states are equilibria. After having done that, we linearize the right-hand side at the splay state and thereby obtain Jacobians with eigenvalues
\begin{align*}
	q_1 = 0, \quad q_{2,3} = \frac K2 e^{\pm i \alpha}
\end{align*}
in a $(1,0)$-phase reduction and
\begin{align*}
	q_1 = 0, \quad q_{2,3} = \frac{K}{2} e^{\pm i\alpha} \left(1 - \frac 1{2m} K e^{\pm i\alpha}\right)
\end{align*}
in a $(2,0)$-phase reduction. An analytical derivation of the eigenvalues in the full system turned out to be too complicated. In both cases, the critical Floquet exponent is given by $q^\mathrm{crit} = q_{2,3}$ and thus the critical PRMM is $\lambda^\mathrm{crit} = e^{T q^\mathrm{crit}}$.

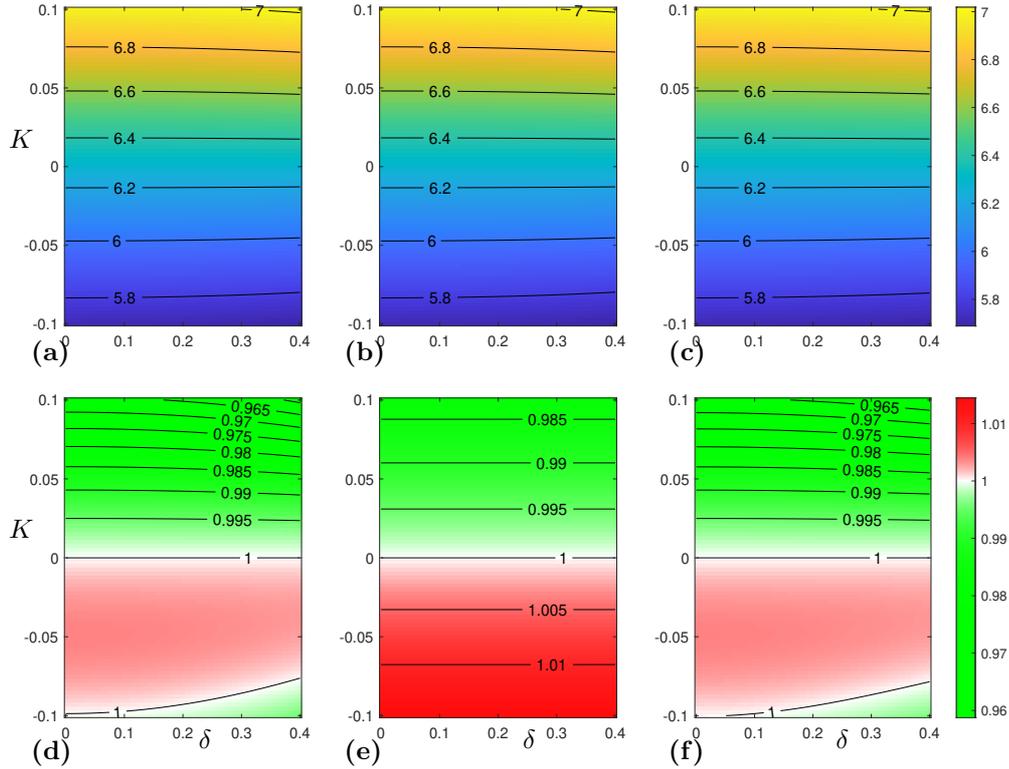
\begin{figure}[ht]
\begin{overpic}[width = \textwidth]{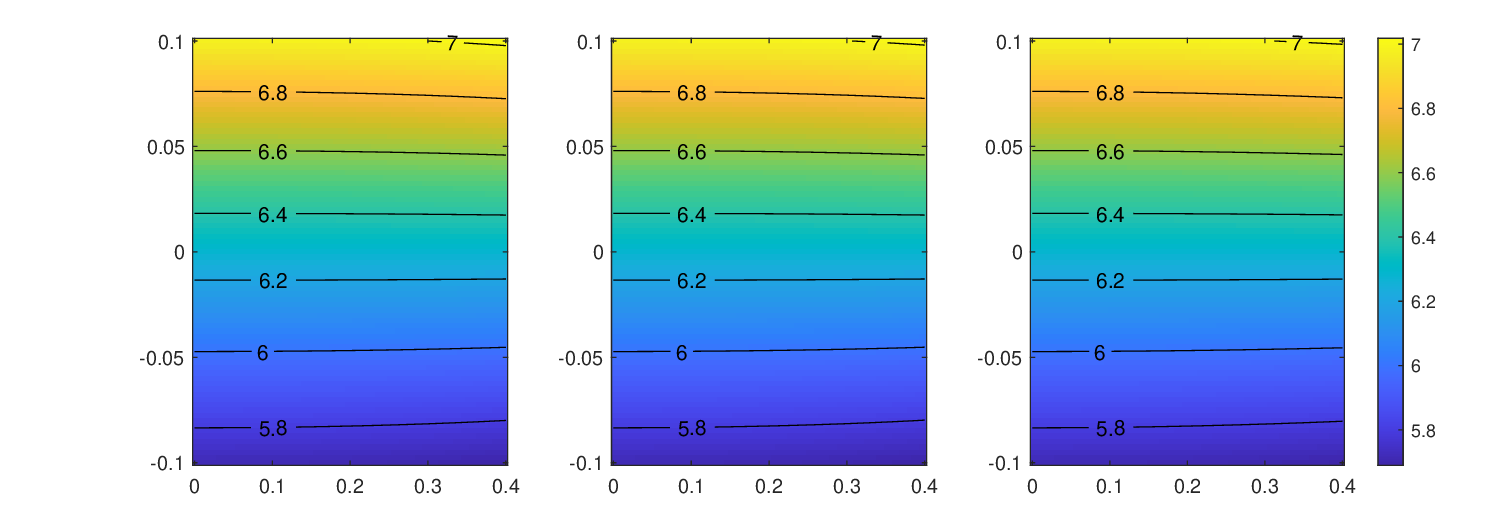}
\put(8,20){$K$}
\end{overpic}

\begin{overpic}[width = \textwidth]{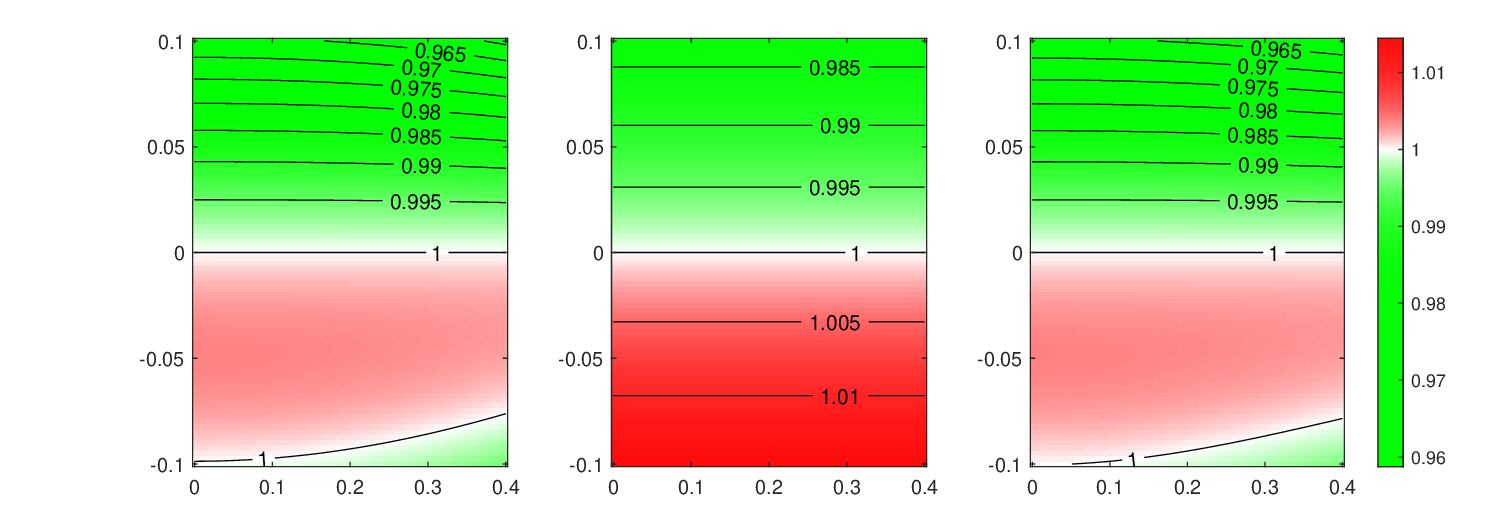}
\put(10,36){\textbf{(a)}}
\put(38,36){\textbf{(b)}}
\put(67,36){\textbf{(c)}}

\put(10,0){\textbf{(d)}}
\put(38,0){\textbf{(e)}}
\put(67,0){\textbf{(f)}}

\put(8,20){$K$}
\put(25,1){$\delta$}
\put(54,1){$\delta$}
\put(82,1){$\delta$}
\end{overpic}

\caption{Numerical calculation of periods and critical multipliers of the Poincar\'e return map of periodic orbits in a neighborhood of the splay orbit. The first column represents the full system, the second column is the $(1,\infty)$-phase reduction and the third column displays the $(2,2)$-phase reduction. The upper row is the period of the resulting periodic orbit, while the lower row depicts the critical PRMM of this orbit. Parameters: $\alpha = \pi/2 +1/20, m = -1, \omega = 1, N=3, g(\phi) = \sin(\phi)$}
\label{fig:floquet_splay}
\end{figure}

While all the previous theory was only valid for $\delta=0$, we now investigate what happens if $\delta \neq 0$. In this case, the right-hand sides of both the phase-reduced systems and the full system are no longer equivariant with respect to~$Z_N$. Therefore, splay states are in general not invariant anymore. However, when $N=3$ and $\delta = 0$, there is a single periodic orbit, i.e., the splay orbit. For general parameter values, this orbit has no Floquet exponents whose absolute value equals~$1$. Thus, this splay orbit is hyperbolic. Slight changes in~$\delta$ away from~$0$ preserve the existence of a periodic orbit in the neighborhood of the splay orbit. To illustrate the stability of these orbits, we numerically search for periodic orbits in a neighborhood of the splay orbit in the full system and phase-reduced systems. Then, we numerically calculate their critical PRMMs to determine the stability of these orbits, see Figure~\ref{fig:floquet_splay}.

A numerical analysis revealed that there is a subcritical Neimark-Sacker bifurcation in Figure~\ref{fig:floquet_splay}(f) when~$K$ is negative and the modulus of the critical PRMM passes through~$1$. In particular, there are two complex conjugated PRMMs that pass through the complex unit circle. This correctly represents the bifurcation behavior of the full system, as in Figure~\ref{fig:floquet_splay}(d). A $(1,\infty)$-phase reduction does not even capture the bifurcation, see Figure~\ref{fig:floquet_splay}(e). The bifurcation at $K=0$ is degenerate. When $K=0$ there is no coupling and so all eigenvalues are $0$.

\section{Phase Reduction Beyond All-To-All Coupled Networks}\label{sec:networks}

In Sections~\ref{sec:meat} and~\ref{sec:dynamics}, we have started with a system of coupled oscillators, derived various phase reductions and compared the stability of synchronized and splay orbits. The system we started with~\eqref{eq:coupling}, consists of $N$ complex Stuart--Landau oscillators coupled with each other via a mean-field coupling, i.e., each oscillator influences every other oscillator in the same way. However, instead of this all-to-all coupling, one could also assume that the coupling between oscillators is described by a graph. In Section~\ref{sec:Networks1} we adopt our calculations from the previous sections to a non all-to-all coupling and Section~\ref{sec:Networks_interpretation} then contains an interpretation of the higher-order interactions terms that we derive from a phase reduction of nonlocally coupled Stuart--Landau oscillators. This provides additional evidence that it is important to study dynamical systems on hypergraphs which have recently been encountered in many different contexts; see for example~\cite{Skardal2019,Bohle2021, Salova2021b,Gong2019,Grilli2017,Gambuzza2021}.

\subsection{Derivation of the phase-reduced Dynamics}\label{sec:Networks1}
 Now, let us assume that the coupling structure is not given by an all-to-all topology, but instead there by a (possibly directed and weighted) graph $\Gamma=(V,E)$ that is described by its adjacency matrix $A\in \R^{N\times N}$ with entries $a_{kl}$. Then, the governing equation is
\begin{align}\label{eq:coupling_graph}
	\dot A_k(t) = \mathcal F(A_k) + K e^{i\alpha}\frac{1}{N}\sum_{l=1}^N a_{kl}(A_l-A_k),
\end{align}
which contains~\eqref{eq:coupling} as a special case when $a_{kl}=1$ for all $k,l=1,\dots,N$. Now, one can do the same analysis as in Section~\ref{sec:meat}, but with~\eqref{eq:coupling_graph} as a starting point. The procedure from Section~\ref{sec:meat} is directly applicable to~\eqref{eq:coupling_graph}, only the resulting formulas are slightly different. Therefore, will not explain the whole procedure again, but only state the results.

The system~\eqref{eq:coupling_graph} can be written in polar coordinates $A_k =r_k e^{i\phi_k}$. Transforming the radii as $R_k = r_k/(1+\delta g(\phi_k))$ yields the system
\begin{subequations}
\label{eq:main_system_graph}
\begin{align}
\label{eq:main_system_R_graph}
\dot R_k &= F(R_k, \phi_k) + KG_k(R,\phi)\\
\label{eq:main_system_phi_graph}
\dot \phi_k &= \omega + K H_k(R, \phi),
\end{align}
\end{subequations}
with functions $F, G_k$ and $H_k$ defined by
\begin{align*}
	F(R_k, \phi_k) &= m R_k^2(R_k-1)(1+\delta g(\phi_k))^2,\\
	G_k(R,\phi) &= \frac{1}{N}\sum_{l=1}^N a_{kl} \Bigg[ R_l\frac{1+\delta g(\phi_l)}{1+\delta g(\phi_k)}\cos(\phi_l-\phi_k+\alpha) - R_k\cos(\alpha)\\
	&\quad - \delta g'(\phi_k)\left( R_l\frac{1+\delta g(\phi_l)}{(1+\delta g(\phi_k))^2}\sin(\phi_l - \phi_k + \alpha) - R_k \frac{\sin(\alpha)}{1+\delta g(\phi_k)}\right) \Bigg],\\
	H_k(R,\phi) &= \frac{1}{N}\sum_{l=1}^N a_{kl}\left[ \frac{R_l(1+\delta g(\phi_l))}{R_k (1+\delta g(\phi_k))}\sin(\phi_l-\phi_k+\alpha) - \sin(\alpha)\right].
\end{align*}
The existence of an invariant torus as in~\eqref{eq:R_expansion} is still guaranteed. Therefore, proceeding as in Section~\ref{sec:meat}, we obtain the first-order phase reduction
\begin{align}
    \nonumber
	\dot \phi_k &=\omega + K P^{(1,\star)}_k(\phi) = \omega + K H_k(1,\phi)\\
    \label{eq:1storder_phase_red_networks}
	&= \omega + K\frac 1N \sum_{l=1}^N a_{kl} \left[ \frac{1+\delta g(\phi_l)}{1+\delta g(\phi_k)}\sin(\phi_l-\phi_k + \alpha) - \sin(\alpha)\right].
\end{align}
To obtain second order phase reductions in~$K$, we need to solve the PDE~\eqref{eq:R1PDE}. When $\delta = 0$, this PDE has the solution
\begin{align*}
	R_k^{(1,0)}(\phi) = \frac{1}{Nm}\sum_{l=1}^N a_{kl} s_0(\phi_k, \phi_l), 
\end{align*}
with $s_0(\phi_k, \phi_l)$ as in~\eqref{eq:R10sol}. Similarly, if $g(\phi)= \sin(\phi)$, on the first order in $\delta$ the solution is
\begin{align*}
	R^{(1,1)}_k(\phi) = \frac{1}{2N(m^2+\omega^2)}\sum_{l=1}^N a_{kl} s_1(\phi_k, \phi_l),
\end{align*}
where $s_1(\phi_k, \phi_l)$ is as in~\eqref{eq:s1}. Finally, if $g(\phi) = \sin(\phi)$, the solution on second order in $\delta$ is
\begin{align*}
	R^{(2,1)}_k(\phi) =	\frac{1}{-4 m N \left(m^4+5 m^2 \omega ^2+4 \omega ^4\right)} \sum_{l=1}^N a_{kl} s_2(\phi_k, \phi_l).
\end{align*}

To determine the second order phase reduction in $K$, one also needs to calculate the gradient of $H$ as illustrated in~\eqref{eq:P2}. It turns out that
\begin{align*}
	\nabla_R H_k^{(-,0)}(1,\phi) = \frac{1}{N}\begin{pmatrix} a_{k1} \sin(\phi_1-\phi_k+\alpha)\\ \vdots \\ a_{kN} \sin(\phi_N-\phi_k+\alpha) \end{pmatrix} - \frac 1N e_k\sum_{l=1}^N a_{kl}\sin(\phi_l-\phi_k+\alpha)
\end{align*}
and that gradients on higher order in $\delta$ are generally of the form
\begin{align*}
	\nabla_R H_k^{(-,\beta)}(1,\phi) = \frac{1}{N}\begin{pmatrix}
	a_{k1} w_\beta(\phi_k, \phi_1)\\ \vdots \\ a_{kN} w_\beta(\phi_k, \phi_N) 	\end{pmatrix}-\frac{1}{N}e_k\sum_{l=1}^N a_{kl} w_\beta(\phi_k, \phi_l),
\end{align*}
for $\beta=0,1,2,\dots$, where $w_\beta$ are trigonometric polynomials. Now, calculating the second order contributions $P_k^{(2,0)}(\phi)$ yields\footnote{The expression for $P_k^{(2,1)}(\phi)$ and $P_k^{(2,2)}(\phi)$ can be generated with the \textsc{Mathematica} code accompanying this paper, see~\cite{Bohle2023GitHub}.}
\begin{subequations}
\label{eq:hypergraph}
\begin{align}
\begin{split}\label{eq:hypergraph1}
	P_k^{(2,0)}(\phi) &= -\frac{1}{2N^2m}\sum_{l=1}^N\sum_{i=1}^N a_{kl}a_{ki}\Big( \sin(\phi_i-\phi_l) - \sin(\phi_i-2\phi_k+\phi_l+2\alpha)\\
	&\qquad\qquad\qquad\qquad\qquad\qquad  - \sin(\phi_k-\phi_l)-\sin(\phi_k-\phi_l-2\alpha) \Big). 
\end{split}\\
\begin{split}\label{eq:hypergraph2}
	&\quad + \frac{1}{2N^2m}\sum_{l=1}^N\sum_{i=1}^N a_{kl} a_{li} \Big( \sin(\phi_i+\phi_k-2\phi_l)-\sin(\phi_k-\phi_l)\\
	&\quad \qquad\qquad\qquad\qquad\qquad-\sin(\phi_k-\phi_l-2\alpha)-\sin(\phi_i-\phi_k+2\alpha) \Big)
\end{split}
\end{align}
\end{subequations}
The next subsection discusses these second-order contributions.

\subsection{Second-Order Phase Reductions as Higher-Order Networks}\label{sec:Networks_interpretation}

We now discuss the individual coupling terms that constitute the $(2,0)$-phase reduction. The coupling includes nonpairwise terms and we discuss how the coupling can be interpreted as a higher-order phase oscillator network on hypergraphs that can be derived from the original graph~$\Gamma = (V,E)$ that describes the coupling of the nonlinear oscillators. In summary, the $(2,0)$-phase reduction is given by
\begin{align*}
    \dot \phi_k = \omega + K P^{(1,0)}_k(\phi) + K^2 P^{(2,0)}_k(\phi),
\end{align*}
where
\begin{align}\label{eq:first-order_graph}
    P^{(1,0)}_k(\phi) = H_k^{(-,0)}(1,\phi) =\frac 1N  \sum_{l=1}^N a_{kl} \Big(\sin(\phi_l-\phi_k+\alpha) - \sin(\alpha)\Big)
\end{align}
agrees with~\eqref{eq:1storder_phase_red_networks} for $\delta = 0$ and~$P^{(2,0)}_k(\phi)$ as specified in~\eqref{eq:hypergraph} contains the second order terms. First, note that the coupling of the first-order phase reduction~\eqref{eq:first-order_graph} is posed on the graph $\Gamma^{(1)} := \Gamma$ that describes the interactions of the coupled nonlinear oscillator network.

\begin{figure}[h]
\centering
\begin{overpic}[width = 0.9\textwidth, grid = false, tics = 5]{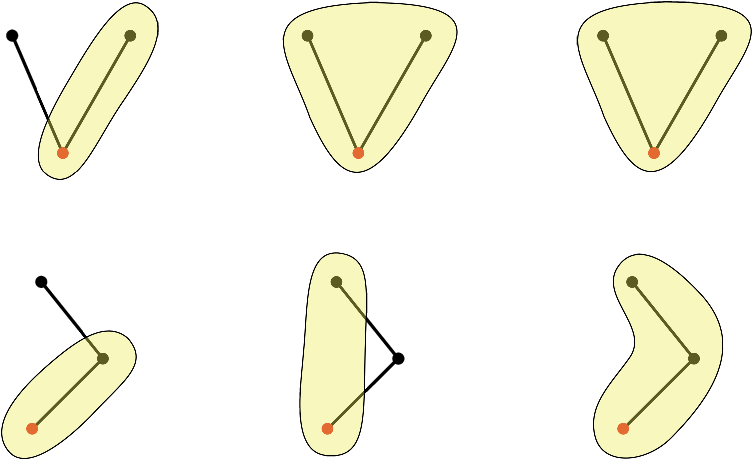}
\put(-2,38){\textbf{(a1)}}
\put(10,41){$k$}
\put(19,56){$l$}
\put(3,56){$i$}
\put(2,62){$\sin(\phi_l-\phi_k+\alpha)$}

\put(-2,-2){\textbf{(a2)}}
\put(6,3){$k$}
\put(15,13){$l$}
\put(7,23){$i$}
\put(2,29){$\sin(\phi_l-\phi_k+\alpha)$}

\put(35,38){\textbf{(b1)}}
\put(49,41){$k$}
\put(58,56){$l$}
\put(42,56){$i$}
\put(42,62){$\sin(\phi_i-\phi_l)$}

\put(35,-2){\textbf{(b2)}}
\put(45,3){$k$}
\put(54,13){$l$}
\put(46,23){$i$}
\put(38, 29){$\sin(\phi_i-\phi_k + 2\alpha)$}

\put(75,38){\textbf{(c1)}}
\put(88,41){$k$}
\put(97,56){$l$}
\put(82,56){$i$}
\put(75,62){$\sin(\phi_i-2\phi_k+\phi_l+2\alpha)$}

\put(75,-2){\textbf{(c2)}}
\put(84,3){$k$}
\put(93,13){$l$}
\put(85,23){$i$}
\put(77,29){$\sin(\phi_i + \phi_k-2\phi_l)$}
\end{overpic}
\caption{Illustration of second order phase interactions in~$K$ that affect node~$k$ (indicated in red).
The first row corresponds to terms appearing in~$\hat g$, while the second row lists terms of~$\bar g$.
The black lines indicate edges in the graph~$\Gamma$ that describe the interaction of the unreduced nonlinear oscillator network. The (hyper)edges that describe the directed phase interactions are indicated by yellow blobs (with node~$k$ being the head).
These include pairwise interactions (panels~a1, a2), pairwise interactions that may be virtual (panel~b2), and three types of nonpairwise interactions that describe the nonlinear influence of nodes~$i,l$ onto node~$k$: One does not depend on the state of node~$k$ itself (panel~b1) and two that depend on the phases of all three nodes (panels~c1, c2).}
\label{fig:HigherOrderInteractions}
\end{figure}
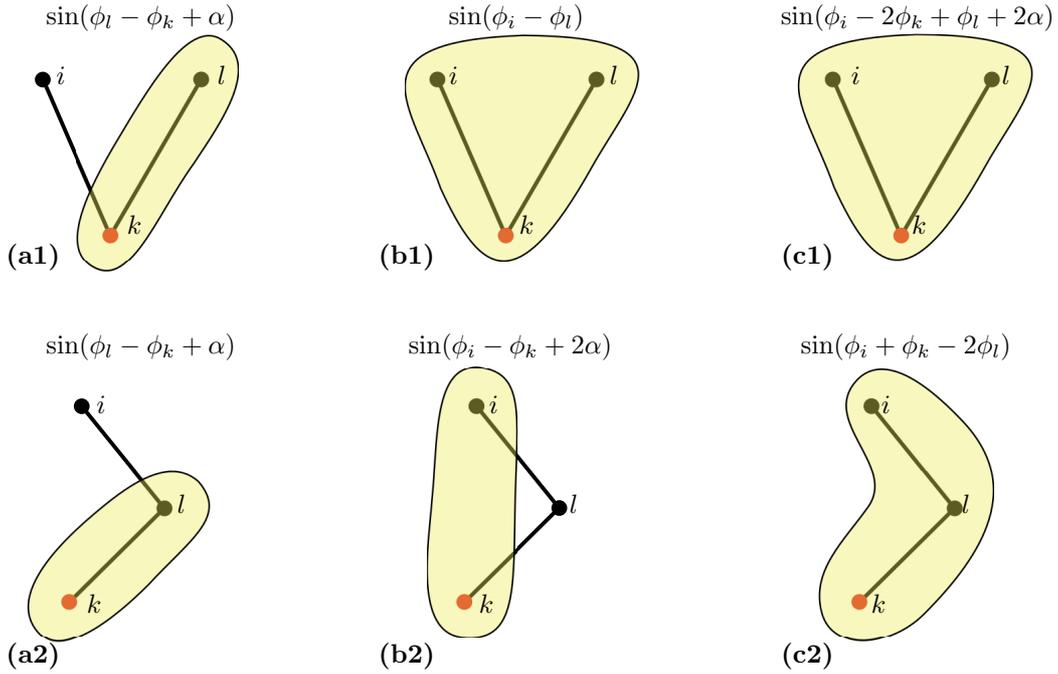

The second-order phase interaction terms~\eqref{eq:hypergraph} not only contain pairwise interactions along network edges but also nonpairwise interactions between triplets of oscillators. 
Collecting all terms, the phase interactions can be interpreted as interactions on two $3$-uniform directed hypergraph with the (adjacency) $3$-tensors $\hat h, \bar h\in \R^{N\times N\times N}$ with coefficients
\begin{align*}
    \hat h_{kli} := a_{kl}a_{ki},\qquad \bar h_{kli} := a_{kl}a_{li},
\end{align*}
where a triplet $(k;l,i)$ corresponds to a directed hyperedge with tail~$k$ and head $\{l,i\}$.

The hypergraph capturing nonpairwise higher-order interactions `inherits' directionality from the original graph~$\Gamma$: The coefficient $h_{kli}$~is nonzero if~$l,i$ are both in the neighborhood of the target node~$k$ in~$\Gamma$ (cf.~Figure~\ref{fig:HigherOrderInteractions}(1)), while $\bar h_{kli}$~is nonzero if there is a path from~$i$ via~$l$ to~$k$ in~$\Gamma$ (cf.~Figure~\ref{fig:HigherOrderInteractions}(2)).
The coupling functions along these hyperedges are
\begin{align}\label{eq:coupling_functions}
\begin{split}
    \hat g(\phi_k, \phi_l, \phi_i) &= 2\cos(\alpha)\sin(\phi_l-\phi_k+\alpha) +\sin(\phi_i-\phi_l) - \sin(\phi_i-2\phi_k+\phi_l+2\alpha),\\
    \bar g(\phi_k, \phi_l, \phi_i) &= 2\cos(\alpha)\sin(\phi_l-\phi_k+\alpha) - \sin(\phi_i-\phi_k+2\alpha) + \sin(\phi_i + \phi_k - 2\phi_l).
\end{split}
\end{align}
These can be obtained from~\eqref{eq:hypergraph} using trigonometric identities and we have
\begin{align*}
    P^{(2,0)}_k(\phi) = -\frac{1}{2N^2m}\sum_{l,i=1}^N \hat h_{kli} \hat g(\phi_k, \phi_l, \phi_i) +\frac{1}{2N^2m}\sum_{l,i=1}^N \bar h_{kli} \bar g(\phi_k, \phi_l, \phi_i).
\end{align*}

Since the coupling functions~\eqref{eq:coupling_functions} contain both pairwise and nonpairwise phase interactions, the interaction structure can be broken down term separating pairwise and nonpairwise interactions: 
Each of the three terms can be associated with a particular type of interaction, which results in six subclasses in total shown in Figure~\ref{fig:HigherOrderInteractions}.

First, there are pairwise correction terms to the first-order phase reduction, that correspond to the first term in the definition of~$\hat g$ and~$\bar g$; see~\eqref{eq:coupling_functions}. 
The coupling of these pairwise correction terms is posed on the graph $\Gamma^{(2)}_\mathrm{a} := \Gamma$  that is the same graph as the coupling of the full system and the coupling of the first-order phase reduction~$\Gamma^{(1)}$; see Figure~\ref{fig:HigherOrderInteractions}(a). 
One of these pairwise correction terms is weighted with the degree of the node~$k$; see Figure~\ref{fig:HigherOrderInteractions}(a1), while the other term is weighted with the degree of node~$l$ (that is in the neighborhood of~$k$), see Figure~\ref{fig:HigherOrderInteractions}(a2). 
Moreover, the interaction function~$\sin(\phi_l-\phi_k+\alpha)$ agrees with the interaction function of a first-order phase reduction~\eqref{eq:first-order_graph} up to the shift~$\sin(\alpha)$ and a constant factor. 
Since the two pairwise correction terms in Figure~\ref{fig:HigherOrderInteractions}(a) share the same coupling structure~$\Gamma^{(2)}_\mathrm{a}$ and the same interaction function, they can be combined into
\begin{align*}
    \frac{-\cos(\alpha) }{N^2m}\left(\sum_{l=1}^N a_{kl} (\deg(k)-\deg(l)) \sin(\phi_l-\phi_k+\alpha)\right).
\end{align*}
Based on this formula, one can see that the sign of the pairwise correction terms is determined by comparing the degree of node~$k$ with the average degrees of all neighbors~$l$ of~$k$.

Second, consider the second term in~$\bar g$ that describes a pairwise interaction between node~$k$ and node~$i$ that is in the neighborhood of a the neighbor~$l$ of~$k$.
While this yields a pairwise interaction from node~$i$ to node~$k$, there may not necessarily be and edge $(i,k)\in E(\Gamma)$ from~$i$ to~$k$ in the graph~$\Gamma$ that describes the original network of coupled nonlinear oscillators.
If $(i,k)\in E(\Gamma)$ then this second-order term describes a second-order correction to the first-order interaction. If $(i,k)\not\in E(\Gamma)$ then this interaction can be considered as a \emph{virtual edge}, which is present in a weighted graph $\Gamma^{(2)}_\mathrm{b2}$ defined by the adjacency matrix $C = (c_{ki})_{k,i=1,\dots,N}$ with coefficients $c_{ki} := \sum_l a_{kl}a_{li}$ but not in~$\Gamma$.
As one can see from the definition of the entries~$c_{ki}$ of the adjacency matrix of~$\Gamma^{(2)}_\mathrm{b2}$, this interaction is weighted by the number of paths connecting~$k$ to~$i$ in the graph~$\Gamma$.

Third, the second term of~$\hat g$ represents a triplet interaction where nodes~$i$ and~$l$---both neighbors of~$k$---influence the node~$k$; see Figure~\ref{fig:HigherOrderInteractions}(b1).
While the interaction is a nonpairwise interaction involving three distinct nodes ($i,l$ and the target~$k$), this may be considered a `nonstandard' nonpairwise interaction as the coupling is independent of the state of node~$k$.
This interaction can be represented by a directed and possibly weighted hypergraph~$\mathcal H^{(2)}_\mathrm{b1}$ represented by a corresponding adjacency $3$-tensor indexed by $k,l,i$.
Note that there will be a symmetry that allows swapping the indices~$l$ and~$i$, but in general this hypergraph is still directed as one can not arbitrarily permute all indices~$k,l,i$.

Finally, there are two further triplet interactions. The first triplet interaction, i.e., the last summand in the definition of~$\hat g$ is an interaction between two neighbors~$i$ and~$l$ of the node~$k$. The coupling structure of this interaction can be described by a directed and weighted hypergraph~$\mathcal H^{(2)}_\mathrm{c1}$, whose $3$-tensor agrees with~$\hat h$. Note that there is a symmetry between~$i$ and~$l$, but one cannot arbitrarily permute all indices which is why in general the hypergraph is directed. The second triplet interaction is governed by the last summand in the definition of~$\bar g$, is one between the node~$k$ itself, a neighbor~$l$ of~$k$ and a neighbor~$i$ of~$l$. Again, the coupling structure can be described by a directed hypergraph~$\mathcal H^{(2)}_\mathrm{c2}$. This time, however, the $3$-tensor that describes the hypergraph corresponds with~$\bar h$ and in general, there this hypergraph does not possess any symmetry with respect to a permutation of indices.

To summarize, the first-order phase reduction, that we have considered, consists of only one interaction term, whose coupling structure~$\Gamma^{(1)}$ agrees with the coupling structure~$\Gamma$ of the full system~\eqref{eq:main_system_graph}. In a second order phase reduction, quite a few new interaction terms appear. While the coupling of some of them is posed on a graph that agrees with~$\Gamma$, the coupling of others is determined by a graph that consists of virtual edges that might not be present in~$\Gamma$. Moreover, there are also three types of triplet interactions on directed hypergraphs, which can be derived from the adjacency matrix of~$\Gamma$.

To conclude this section, we want to remark that a second order phase reduction contains interaction terms on directed hypergraphs, even when the underlying graph~$\Gamma$, that determined the coupling in the full system~\eqref{eq:main_system_graph}, is undirected and unweighted. The only exception is when~$\Gamma$ itself is an all-to-all graph. However, whenever~$\Gamma$ is connected, yet non all-to-all, there exists an open triangle as seen in Figure~\ref{fig:HigherOrderInteractions}(c1), which causes the hypergraph, that governs the second order phase reduction, to be directed.

\section{Discussion}\label{sec:conclusion}

Phase reductions provide a useful tool to analyze the dynamics of coupled oscillator networks. 
Here we derived explicit expressions for nonlinear oscillations with phase-dependent amplitude subject to simple diffusive coupling.
By using a suitable coordinate transformation, our results also apply to systems where the limit cycle is simple but the coupling is phase dependent or a combination thereof.
The reduced phase equations allow to analyze, for example, phase dynamics for oscillators that are---if uncoupled---further away from a Hopf bifurcation point where the limit cycle becomes noncircular in general.

While the shape of the limit cycle affects the collective dynamics, a first-order phase reduction is insufficient to capture the dynamical effects of the amplitude dependence:
The phase reduction needs to be at least of second order in both the coupling strength~$K$ and the parameter~$\delta$ that describes the perturbation from a circular limit cycle.
We showed that second-order phase reductions were able to accurately predict the stability properties of the synchronized and splay orbit when all terms of up to second order in~$K$ and~$\delta$ were included.
Importantly, the amplitude dependence breaks the rotational symmetry of phase equations that is typical, for example, of the Kuramoto equations. 
While such symmetry breaking has been analyzed from the perspective of the phase equations~\cite{Brown2003}, in our setup it arises through a perturbation of the underlying nonlinear oscillator.
This perspective allows to make direct comparisons between the nonlinear system and the phase-reduced dynamics in contrast to phase reductions via normal forms~\cite{Ashwin2016a}, where normal form symmetries---that may be absent in the full equations~\cite{Crawford1991}---can appear in the phase reduction.

As detailed in Section~\ref{sec:Networks_interpretation}, the second-order phase interaction term for coupled oscillators on a given graph~$\Gamma$ can be interpreted as phase oscillator dynamics on (directed) hypergraphs. 
Specifically, the second-order phase interaction terms correspond to second-order corrections to interactions along edges of~$\Gamma$, possible virtual pairwise connections between oscillator pairs that are not joined by an edge in~$\Gamma$, and nonpairwise triplet interactions of different type.
This highlights the importance to consider the dynamics on directed hypergraphs---these may appear implicitly (e.g.,~\cite{Bick2017c}) while explicit frameworks with distinct perspectives have been introduced in~\cite{Aguiar2020,Gallo2022,VonderGracht2023b}.
Indeed, directionality of higher-order interactions has implications for synchrony~\cite{Gallo2022} or the emergence of more complicated dynamical phenomena such as heteroclinic dynamics~\cite{Bick2017c,Bick2023b}.

In principle, the analysis presented in Section~\ref{sec:meat} can be extended to derive higher-order interactions beyond second order in both the coupling strength~$K$ and the deviation parameter~$\delta$.
Such higher-order interactions would include non-additive interactions between quadruplets of phases and allow to describe the approximate phase dynamics beyond second order.
The main obstacle that has to be overcome is the algebraic complexity of these phase reductions. 
Already at second order, the terms for triplet interactions become quite complex that bring symbolic computer algebra software, such as \textsc{Mathematica}, to their limits.

Generalizing the computations to oscillator networks with nonlinear coupling, nonidentical oscillators, or strongly nonsinusoidal oscillations---all properties of real-world oscillator networks---highlights some of the challenges to compute phase reductions explicitly.
First, nonlinear coupling makes equation~\eqref{eq:higher-order_expans} more challenging to solve (if an explicit solution is possible at all).
If the nonlinear coupling contains nonlinear coupling between three oscillators (e.g., a coupling terms of the form~$A_lA_jA_k$)---corresponding to then we expect nonpairwise coupling already in the first-order phase reduction (see~\cite{Ashwin2016}).
If the coupling is nonlinear in one oscillator, i.e., the coupling only includes coupling terms of the type~$A_j^m$ (as, for example, in~\cite{Kori2008}) then we have higher harmonics in~\eqref{eq:higher-order_expans} making it less tractable.
Second, real-world oscillators are rarely perfectly identical motivating the question how heterogeneity affects the phase reduction.
Here we assumed that the oscillators are identical and, in particular, that~$\omega$ in~\eqref{eq:main_system_phi} is independent of~$k$. 
If~$\omega$ depended on~$k$ the first-order phase reduction would not change at all and one could just replace~$\omega$ by~$\omega_k$ everywhere. 
A second-order phase reduction, could theoretically also be derived by adapting the methods from Section~\ref{sec:meat}. 
Practically, however, the terms of second order start to depend nonlinearly on $\omega=(\omega_1,\dots,\omega_N)$. Thus, already when $\delta=0$ and~$N$ is small and finding a general solution for $R_k^{(1,0)}(\phi)$ is challenging. 
Moreover, due to the growing complexity of the second-order phase reduction, it is practically intractable. 
Extending this to $\delta >0$ only worsens the problem. A possible approach to overcome this problem is to assume that the intrinsic frequencies~$\omega_k$ are sampled from a probability distribution and consider a mean-field limit. 
Third, strongly nonlinear oscillations, such as FitzHugh--Nagumo neurons, are far away from a circular limit cycle (and would this require a large perturbation parameter~$\delta$). While this suggests the importance of expanding to high-order in~$\delta$, more direct approaches may be more appropriate~\cite{Izhikevich2000} that allow to get a more qualitative understanding of the dynamics~\cite{Ashwin2021}.

\paragraph{Acknowledgements}
We thank P.~Ashwin, H.~Nakao, and  M.~T.~Schaub for helpful discussions.
We gratefully acknowledge the support of the Institute for Advanced Study at the Technical University of Munich through a Hans Fischer Fellowship awarded to CB that made this work possible. CB acknowledges support from the Engineering and Physical Sciences Research Council (EPSRC) through the grant EP/T013613/1. TB acknowledges support of the TUM TopMath elite study program. CK acknowledges support via a Lichtenberg Professorship of the VolkswagenStiftung.

\paragraph{Data Availability}
The \textsc{Matlab} code that generates the figures and the \textsc{Mathematica} code that computes the phase reductions is publicly available on a GitHub repository that can be accessed via~\href{https://github.com/tobiasboehle/HigherOrderPhaseReductions}{https://github.com/tobiasboehle/HigherOrderPhaseReductions}, Ref.~\cite{Bohle2023GitHub}

\pagebreak

\appendix

\section{$(2,1)$-Phase Reduction for Arbitrary Perturbations $g$}

In this section, we state the solution $R_k^{(1,1)}$ of the PDE~\eqref{eq:PDE_d1} when $g$ is not just given by $\sin$ but consists of more harmonics. We explain how higher-harmonics in $g$ influence the $(2,1)$-phase reduction and how these additional harmonics affect the stability of the synchronized orbit.

Whenever $g(\phi) = \sin(n \phi)$ for $n\in\N$, the solution of the PDE~\eqref{eq:PDE_d1} is given by
\begin{align}\label{eq:R11solgen}
	R^{(1,1)}_k(\phi) = \frac{1}{2N(m^2 + (n\omega)^2)}\sum_{l=1}^N s_1(\phi_k, \phi_l)
\end{align}
with
\begin{align*}
	s_1(\phi_k, \phi_l) &= n\omega\Big((n-2)\cos(n\phi_k-\alpha) - \cos(\phi_k-(n+1)\phi_l-\alpha)\\
	&\qquad- (n-3)\cos((n+1)\phi_k-\phi_l-\alpha) - \cos(\phi_k + (n-1)\phi_l-\alpha)\\
	&\qquad- (n+2)\cos(n\phi_k+\alpha) + (n+3)\cos((n-1)\phi_k+\phi_l+\alpha)\Big)\\
	&+ m\Big((n-2)\sin(n\phi_k-\alpha) + \sin(\phi_k - (n+1)\phi_l - \alpha)\\
	&\qquad- (n-3)\sin((n+1)\phi_k-\phi_l-\alpha) - \sin(\phi_k+(n-1)\phi_l-\alpha)\\
	&\qquad- (n+2)\sin(n\phi_k+\alpha) + (n+3)\sin((n-1)\phi_k+\phi_l+\alpha)\Big).
\end{align*}
Moreover, if $g(\phi) = \cos(n\phi)$ for $n\in \N$, the solution of~\eqref{eq:PDE_d1} is given by~\eqref{eq:R11solgen} as well, but then
\begin{align*}
	s_1(\phi_k, \phi_l) &= n\omega \Big(-(n-2)\sin(n\phi_k-\alpha) - \sin(\phi_k-(n+1)\phi_l-\alpha)\\
	&\qquad + (n-3)\sin((n+1)\phi_k - \phi_l-\alpha) + \sin(\phi_k+(n-1)\phi_l-\alpha)\\
	&\qquad+ (n+2)\sin(n\phi_k+\alpha) - (n+3)\sin((n-1)\phi_k+\phi_l+\alpha) \Big)\\
	&+ m\Big((n-2)\cos(n\phi_k-\alpha) - \cos(\phi_k-(n+1)\phi_l-\alpha)\\
	&\qquad-(n-3)\cos((n+1)\phi_k-\phi_l-\alpha) - \cos(\phi_k + (n-1)\phi_l-\alpha)\\
	&\qquad- (n+2)\cos(n\phi_k + \alpha) + (n+3)\cos((n-1)\phi_k+\phi_l+\alpha)\Big)
\end{align*}

Now, a general sufficiently smooth function $g$ can be constructed as a sum of the basis functions $\cos(n\phi)$ and $\sin(n\phi)$ with $n\in\N$. Due to the linearity of the PDE~\eqref{eq:PDE_d1} its solution for a general sufficiently smooth function $g$ can therefore be constructed using its solutions when~$g$ is a basis function.

To investigate how these higher harmonics influence the stability of the synchronized orbit in a $(2,1)$-phase-reduced system one first notes that the only part in this phase reduction, that depends on $R^{(1,1)}_k(\phi)$ is $P^{(2,1)}_k$ as defined in~\eqref{eq:P21}. In particular, if $g$ is given by a general Fourier sum
\begin{align*}
	g(\phi) = \sum_{n=1}^M (a_n \cos(n\phi) + b_n \sin(n\phi))
\end{align*}
with real coefficients $a_n, b_n$, the general form of $P_k^{(2,1)}$ is
\begin{align*}
	P_k^{(2,1)}(\phi) = \nabla_R H_k^{(-,0)}(1,\phi)\cdot \left(\sum_{n=1}^M a_n R^{(1,1)}_{\cos(n\phi)}(\phi) + b_n R^{(1,1)}_{\sin(n\phi)}(\phi)\right) + \nabla_R H_k^{(-,1)}(1,\phi)\cdot R^{(1,0)}(\phi),
\end{align*}
where $R^{(1,1)}_{\cos(n\phi)} = (R^{(1,1)}_{\cos(n\phi),1},\dots,R^{(1,1)}_{\cos(n\phi),N})^\top$ and $R^{(1,1)}_{\sin(n\phi)} = (R^{(1,1)}_{\sin(n\phi),1},\dots, R^{(1,1)}_{\sin(n\phi),N})^\top$ are the solutions~\eqref{eq:R11solgen} of~\eqref{eq:PDE_d1} when $g(\phi) = \cos(n\phi)$ and $g(\phi) = \sin(n\phi)$, respectively. In other words, one can say that whenever $g$ consists of multiple harmonics, these harmonics contribute to the right-hand side of the $(2,1)$-phase reduction, each by one summand $\nabla_R H_k^{(-,0)}(1,\phi)\cdot R^{(1,1)}_{\cos(n\phi)}(\phi)$ or $\nabla_R H_k^{(-,0)}(1,\phi)\cdot R^{(1,1)}_{\sin(n\phi)}(\phi)$ with possible prefactors, only. Therefore, higher harmonics in~$g$ cause more summands in the linearization of the right-hand side at a synchronized state. As explained in Section~\ref{sec:dynamics_sync} this Jacobian is of the form $h(\gamma) \frac{1}{N}(\mathbbm 1 - N \mathbb I)$, when linearized at $\phi_1=\dots=\phi_N=\gamma$ and each harmonic in~$g$ causes one summand in $h$. A calculation shows that this summand for the harmonic $\sin(n\phi)$ is 
\begin{align*}
	\frac{2 K^2 \delta \sin(\alpha)^2}{m^2 + (n\omega)^2} \Big( n\omega \cos(n\gamma) + m \sin(n\gamma)\Big)
\end{align*}
and
\begin{align*}
	\frac{2 K^2 \delta \sin(\alpha)^2}{m^2 + (n\omega)^2} \Big( m\cos(n\gamma) - n\omega \sin(n\gamma)\Big)
\end{align*}
if $g(\phi) = \cos(n\phi)$. When integrating these summands over the synchronized orbit $\gamma\in \S$ one sees that they vanish. Therefore, they do not contribute to the Floquet exponent, which determines the stability of the synchronized orbit.

\bibliographystyle{plain}
\bibliography{library}

\end{document}